\documentclass[twoside]{amsart}
	\usepackage{amsfonts,amsmath,amssymb,amsthm}
	\usepackage{mathrsfs,enumerate,color,datetime,pxfonts,extarrows,graphicx,subcaption}
	\usepackage[all]{xy}
	\usepackage[utf8]{inputenc}
	\usepackage[colorlinks,citecolor=blue]{hyperref}

	\renewcommand{\AA}{\mathbf{A}}
	\newcommand{\BB}{\mathbf{B}}	
	\newcommand{\CC}{\mathbf{C}}	
	\newcommand{\FF}{\mathbf{F}}
		
	\newcommand{\NN}{\mathbf{N}}	
	\newcommand{\PP}{\mathbf{P}}
	\newcommand{\RR}{\mathbf{R}}
	\newcommand{\QQ}{\mathbf{Q}}	
	\newcommand{\WW}{\mathbf{W}}	
	\newcommand{\ZZ}{\mathbf{Z}}	
	\newcommand{\Cc}{\mathscr{C}}
	\newcommand{\Ff}{\mathscr{F}}	
		
	\newcommand{\Hh}{\mathscr{H}}	
	\newcommand{\Ll}{\mathscr{L}}	
	\newcommand{\Oo}{\mathscr{O}}	
	\newcommand{\Xx}{\mathscr{X}}	
	\newcommand{\mm}{\mathfrak{m}}	

	\newcommand{\et}{{\rm\acute{e}t}}
	\newcommand{\rig}{{\rm rig}}
	\newcommand{\id}{{\rm id}}
	\DeclareMathOperator{\Spec}{Spec}
	
	\DeclareMathOperator{\Spf}{Spf}
	\DeclareMathOperator{\Spa}{Spa}
	\DeclareMathOperator{\Swan}{Sw}
	\DeclareMathOperator{\dimtot}{dimtot}	
	
	\DeclareMathOperator{\fil}{fil}
	\DeclareMathOperator{\Gr}{Gr}
	\DeclareMathOperator{\rk}{rk}
	\DeclareMathOperator{\Hom}{Hom}
	
	\newcommand{\ANK}[2]{{\AA_{#2}}^{\! \! \! \!  #1}} 
	\newcommand{\AN}[1]{\ANK{#1}{k}} 

	\newcommand{\isom}{\simeq}
	\newcommand{\xto}{\xlongrightarrow}
	\renewcommand{\phi}{\varphi}
	\renewcommand{\to}{\longrightarrow}
	\renewcommand{\bar}{\overline}
	\newcommand*\isomlong{%
 		\xlongrightarrow{\raisebox{-0.2 em}{\smash{\ensuremath{\sim}}}}%
	}

	\theoremstyle{plain}
		\newtheorem{lemma}[subsubsection]{Lemma}
		\newtheorem{proposition}[subsubsection]{Proposition}
		\newtheorem{theorem}[subsubsection]{Theorem} 
		\newtheorem{corollary}[subsubsection]{Corollary}
		\newtheorem{question}[subsubsection]{Question}
		\newtheorem*{proposition*}{Proposition}
		\newtheorem*{theorem*}{Theorem}
	\theoremstyle{definition}
		\newtheorem{definition}[subsubsection]{Definition}
	\theoremstyle{remark}
		\newtheorem{example}[subsubsection]{Example}
		\newtheorem{remark}[subsubsection]{Remark}

	\numberwithin{equation}{subsection}

	\author[P. Achinger]{Piotr Achinger}
	\address{Instytut Matematyczny PAN, ul. Śniadeckich 8, 00-656 Warszawa, Poland}  
	\address{Institut des Hautes \'Etudes Scientifiques \\
		Le Bois-Marie 35, route de Chartres \\
		91440 Bures-sur-Yvette, France}
	\email{achinger@ihes.fr}
	\title{Wild ramification and $K(\pi, 1)$ spaces}
	\date{\today}
\begin{document}

\begin{abstract}
	We prove that every connected affine scheme of positive characteristic is a $K(\pi, 1)$ space for the \'etale topology. The main ingredient is the special case of the affine space $\AN n$ over a field $k$. This is dealt with by induction on $n$, using a key ``Bertini-type'' statement regarding the wild ramification of $\ell$-adic local systems on affine spaces, which might be of independent interest. Its proof uses in an essential way recent advances in higher ramification theory due to T.~Saito. We also give rigid analytic and mixed characteristic versions of the main result.
\end{abstract}

\maketitle

\section{Introduction} \label{s:intro}

\stepcounter{subsection}

The \'etale homotopy theory of schemes of positive characteristic is quite poorly understood. For example, already the \'etale fundamental group $\pi_1(\AN 1)$ of the affine line over an algebraically closed field $k$ of characteristic $p$ is incredibly complicated \cite{Raynaud}.  One of our main results is the realization that the \'etale homotopy theory of characteristic $p$ schemes is in a certain way controlled by the \'etale fundamental group. 

The notion of a $K(\pi, 1)$ space for the \'etale topology plays a central role in this paper. In algebraic topology, a connected topological space $X$ (with a basepoint $x$) satisfying suitable technical assumptions is called a \emph{$K(\pi, 1)$ space} if its higher homotopy groups $\pi_q(X,x)$ ($q>1$) are zero. The homotopy type of such a space is completely determined by its fundamental group $\pi_1(X,x)$, and in particular the cohomology of any local system on $X$ agrees with the group cohomology of the corresponding representation of $\pi_1(X,x)$. Similarly, we call a connected scheme $X$ with a geometric point $\bar x$ a \emph{$K(\pi, 1)$ scheme} if for every locally constant \'etale sheaf of finite abelian groups $\Ff$ on $X$, the natural maps
\begin{equation} \label{eqn:rhomaps-intro} 
	H^*(\pi_1(X, \bar x), \Ff_{\bar x}) \to H^*(X, \Ff)
\end{equation}
are isomorphisms (cf. Definition~\ref{def:kpi1}). If $X$ is noetherian and geometrically unibranch, this is equivalent to the vanishing of the higher \'etale homotopy groups $\pi_q(X, \bar x)$ for $q>1$ (cf. Proposition~\ref{prop:kpi1artinmazur}).

One of our main results is the following.

\begin{theorem} \label{thm:affine-kpi1}
	Every connected affine $\FF_p$-scheme is a $K(\pi, 1)$ scheme.
\end{theorem}

This of course stands in stark contrast with the characteristic zero case. One might interpret this by saying that the \'etale fundamental group of a connected affine $\FF_p$-scheme is so large that it `absorbs' the higher homotopy groups. To go one step beyond the affine case, Theorem~\ref{thm:affine-kpi1} implies that the \'etale homotopy type of a normal quasi-compact and separated $\FF_p$-scheme can be described as the homotopy colimit of a finite diagram of classifying spaces of profinite groups.

\subsection{Artin neighborhoods}
\label{ss:artin}

To explain both our initial motivation and the idea of proof of Theorem~\ref{thm:affine-kpi1}, we will start by discussing Artin's construction of $K(\pi, 1)$ neighborhoods on smooth complex algebraic varieties. 

The notion of a $K(\pi, 1)$ scheme first appeared in algebraic geometry in Artin's proof of the comparison theorem between singular and \'etale cohomology of a smooth scheme $X$ over $\CC$, which states that
\begin{equation} \label{eqn:etale-singular}
	H^*(X, \Ff) \isomlong H^*(X(\CC), \Ff)
\end{equation}
for every locally constant \'etale sheaf of finite abelian groups $\Ff$ on $X$ \cite[Exp.~XI]{Artin}.

Recall that an \emph{elementary fibration} is a morphism of schemes $\pi\colon X\to S$ for which there exists a commutative diagram
\[ 
	\xymatrix{
		X \ar[dr]_\pi \ar[r]^j & \bar X \ar[d]^{\bar \pi} & Y \ar[dl] \ar[l]_i \\
		& S &
	}
\] 
where $j\colon X\to \bar X$ is an open immersion, $i\colon Y\to \bar X$ the complementary closed immersion, $\bar \pi$ is projective with geometrically connected fibers, smooth of dimension $1$, and its restriction $\bar \pi|_Y \colon Y\to S$ is finite \'etale surjective. Thus the geometric fibers of $\pi$ are smooth affine curves. Moreover, if $S$ is of finite type over $\CC$, then the associated morphism $\pi\colon X(\CC)\to S(\CC)$ is a locally trivial fibration. Artin showed that a smooth variety $X$ over an infinite field $K$ can be covered by Zariski open subsets $U$ for which there exists a chain of elementary fibrations
\[ 
	U = U_{\dim X} \to \ldots \to U_1 \to U_0 = \Spec K.
\]
We call a scheme $U$ admitting a chain as above an \emph{Artin neighborhood}. Thus $U$ is an `iterated fibration in affine curves,' and it follows that if $K$ has characteristic zero, an Artin neighborhood is a $K(\pi, 1)$ scheme. Moreover, if $K=\CC$, then the associated analytic space $U(\CC)$ is a $K(\pi, 1)$ space, and $\pi_1(U(\CC))$ is an iterated extension of free groups. Consequently, in the natural commutative square
\[ 
	\xymatrix{
		H^*(\pi_1(U,u), \Ff_{u}) \ar[r] \ar[d] & H^*(\pi_1(U(\CC), u), \Ff_u) \ar[d] \\
		H^*(U, \Ff) \ar[r] & H^*(U(\CC), \Ff)
	}
\]
the vertical and top arrows are isomorphisms, which implies \eqref{eqn:etale-singular} for $X=U$ an Artin neighborhood. The general case of \eqref{eqn:etale-singular} follows then easily by cohomological descent.

Our initial goal in this project was to generalize Artin's theorem by showing that a smooth scheme over an infinite field of positive characteristic admits a covering by $K(\pi, 1)$ open subschemes; we dared not hope that something as striking as Theorem~\ref{thm:affine-kpi1} can be true. To this end, a good understanding of the problem with extending Artin's characteristic zero proof put us on the right track.

Namely, the reason why the argument that an Artin neighborhood is a $K(\pi, 1)$ scheme works only in characteristic zero is related to wild ramification. Suppose that $\pi\colon X\to S$ is an elementary fibration, and we know that $S$ is a $K(\pi, 1)$. To show that $X$ has to be a $K(\pi,1)$ as well, it is easily seen that it is enough to prove that for a locally constant constructible sheaf $\Ff$ on $X$, the following condition is satisfied:
\begin{equation} \label{eqn:pi-good-for-F}
	\begin{aligned}
		&\text{The higher direct images } R^q \pi_* \Ff\text{ are locally constant,} \\ 
		&\text{with formation commuting with base change for all }q\geq 0. 
	\end{aligned}
\end{equation}
This last assertion is true in characteristic zero e.g. by comparison with the complex case, where $\pi$ is topologically a locally trivial fibration, but in positive characteristic this can fail due to phenomena related to the wild ramification of the restriction of $\Ff$ to the fibers of $f$ at the points of $Y$, as the following basic example shows.

\begin{example} \label{ex:basic}
Let ${\rm pr}\colon \AN 2 \to \AN 1$ be the projection to the second coordinate, where $k$ is algebraically closed of characteristic $p>0$, $\Ff$ the rank one sheaf of $\FF_\ell$-vector spaces ($\ell\neq p$) on $\AN 2$ associated to the Artin--Schreier covering
\[ 
	t^p - t = x_1 x_2
\]
and a nontrivial character $\psi\colon \FF_p \to \FF_\ell^\times$. Then the restriction of $\Ff$ to $\pi^{-1}(0)$ is constant, while for $x_2\neq 0$, $\Ff|_{\pi^{-1}(x_2)}$ is not constant. Consequently, $\pi_* \Ff = 0$ does not have the required compatibility with base change, and \eqref{eqn:pi-good-for-F} does not hold.
\end{example}

\subsection{The Bertini theorem for lcc sheaves}

Our main technical result below states that one can make problems as in Example~\ref{ex:basic} go away by applying a non-linear automorphism of the affine space. Consequently, one can make the above inductive argument work for sheaves on $\AN {n+1}$ if one is allowed to choose the fibration $\pi\colon \AN{n+1} \to \AN n$ \emph{after} being given the locally constant sheaf $\Ff$. For brevity, let us call a sheaf $\Ff$ \emph{well-aligned} with respect to a map $\pi$ if condition \eqref{eqn:pi-good-for-F} holds.

\begin{theorem}[Bertini theorem for lcc sheaves] \label{thm:bertini}
	Let $k$ be an infinite field of characteristic $p>0$. Let $\ell \neq p$ be a prime and let $\Ff$ be a locally constant constructible sheaf of $\FF_\ell$-vector spaces on $\AN{n+1}$ ($n\geq 0$). Let $\pi\colon \AN{n+1}\to \AN{n}$ be the projection to the first $n$ coordinates. Then there exists an automorphism $\phi$ of $\AN{n+1}$ such that $\phi^* \Ff$ is well-aligned with respect to $\pi$ (equivalently, $\Ff$ is well-aligned with respect to $\pi\circ \phi^{-1}$). 
\end{theorem}

In fact, we prove that the assertion holds for a (not necessarily linear!) automorphism $\phi$ which is general in a suitable sense. This justifies the name `Bertini theorem.' See Theorem~\ref{thm:bertini-general} for a precise formulation. 

We shall now explain how Theorem~\ref{thm:bertini} implies Theorem~\ref{thm:affine-kpi1}. First, Theorem~\ref{thm:bertini} enables us to prove by induction on $n$ that the affine space $\AN n$ is a $K(\pi, 1)$, along the lines sketched in \S\ref{ss:artin}. This turns out to be the key case of Theorem~\ref{thm:affine-kpi1}. To deduce the general case, one first treats affine \'etale schemes over $\AN n$ as an intermediate step. To this end, one uses the following result.

\begin{proposition*}[\ref{prop:nagatatrick}]
	Let $U$ be an affine scheme of finite type over $k$ admitting an \'etale map $g\colon U\to \AN{n}$. Then there exists a finite \'etale map $f\colon U\to \AN{n}$.
\end{proposition*}

The proof of this assertion is suprisingly easy, and is based on Nagata's proof of the Noether normalization lemma \cite[I \S 1]{Mumford}. A similar result has been obtained by Kedlaya \cite{Kedlaya}, and both were inspired by a trick used in \cite{Katz} in the one-dimensional case. The following example illustrates the general idea.

\begin{example}
 	Let $U = \AN 1 \setminus \{x_1, \ldots, x_r\}$ and let $g\colon U\to \AN 1$ be the inclusion. Then $g$ is of course \'etale, but it is not finite because the points $x_i$ have empty preimages. To remedy this, we can send these points off to infinity by adding to $g$ functions with poles at the $x_i$. This will make the map $g$ finite, but might destroy \'etaleness --- unless we are in characteristic $p$, in which case the added functions can be taken to be $p$-th powers. Concretely, we might take
 	\[ 
 		f(x) = x + \sum_{i=1}^r \frac{1}{(x-x_i)^p}.
 	\]
\end{example} 

One can apply a similar reasoning in a mixed characteristic situation and prove that an affinoid rigid space which is \'etale over a polydisc is also finite \'etale over a polydisc, cf. Proposition~\ref{prop:nagata-rigid}.

Combining Proposition~\ref{prop:nagatatrick} with the fact that $\AN n$ is a $K(\pi, 1)$, we see that if $U$ is an affine scheme of finite type over $k$, admitting an \'etale map $U\to \AN n$, then $U$ is a $K(\pi, 1)$. Finally, using limit arguments and Gabber's affine analog of the proper base change theorem \cite{Gabber}, we deduce Theorem~\ref{thm:affine-kpi1}.

\subsection{Higher ramification theory and proof of the Bertini theorem}

Theorem~\ref{thm:bertini} is where higher ramification theory enters the picture. A key ingredient in the proof is the Deligne--Laumon theorem (cf. Corollary~\ref{cor:deligne-laumon}), which yields a condition for the higher direct images $R^q \pi_* \Ff$ ($q\geq 0$) being locally constant in terms of the Swan conductors at infinity of the restriction of $\Ff$ to the fibers of $\pi$. The `baby case' is when $\Ff$ is non-fierce at infinity (cf. Definition~\ref{def:nonfierce-sheaf}), in which case we can take $\phi$ to be a general linear automorphism (cf. Proposition~\ref{prop:bertini-nonfierce}). In the general case, we use the recent work of Takeshi Saito on the characteristic cycle associated to a locally constant $\FF_\ell$-sheaf \cite{Saito,Saito2016}. It turns out that we can take the automorphism $\phi$ to be quadratic. 

As pointed out to us by Maxim Kontsevich, it makes sense to ask whether a variant of Theorem~\ref{thm:bertini} holds for irregular connections on $\ANK{n+1}{\CC}$. Our method of proof, employing the characteristic cycle, seems to suggest that the answer should be positive. We plan to address this question in a future paper.

In the course of our work on Theorem~\ref{thm:bertini}, we started by solving its rank one case first. In this situation, the calculations are very explicit, and we include them in an appendix. One good feature of this proof is that, unlike our treatment of the general case, the arguments work over a finite field. It could be interesting to obtain a general proof of our Bertini theorem over finite fields.

\subsection{Mixed characteristic and rigid analytic variants}

If $X$ is a connected affine $\FF_p$-scheme and $\Ff$ a locally constant constructible $\FF_p$-sheaf on $X$, then the maps \eqref{eqn:rhomaps-intro} are isomorphisms. This can be seen easily using the Artin--Schreier sequence, cf. Example~\ref{ex:padickpi1}. Scholze \cite[Theorem~4.9]{ScholzePAdic} observed that using perfectoid spaces, one can deduce that every mixed characteristic noetherian affinoid adic space is a $K(\pi, 1)$ space for $p$-adic coefficients. Using a similar argument and Theorem~\ref{thm:affine-kpi1}, one can give the following strenghtening of Scholze's result.

\begin{theorem*}[\ref{thm:kpi1affinoid-mixed}]
	Every noetherian affinoid adic space over ${\rm Spa}(\QQ_p, \ZZ_p)$ is a $K(\pi, 1)$ space. 
\end{theorem*}

This in turn allows us to give a mixed characteristic variant of Theorem~\ref{thm:affine-kpi1}.

\begin{theorem*}[\ref{thm:henselian-mixed}]
	Let $A$ be a noetherian $\ZZ_{(p)}$-algebra such that $(A, pA)$ is a henselian pair. Then $\Spec A$ and $\Spec A[1/p]$ are $K(\pi, 1)$ schemes.
\end{theorem*}

This has a natural application to Milnor fibers and Faltings' topos, allowing us to remove the log smoothness hypothesis of the main result of \cite{Achinger}, cf. Corollary~\ref{cor:better-achinger}.

\subsection{The naive \'etale topology}

Let $X$ be a scheme. By the \emph{naive \'etale topology} we mean the topology on the category of \'etale $X$-schemes generated by Zariski coverings and finite \'etale surjective maps. The corresponding topos, which we denote by $X_{{\rm n.}\et}$, is related to the \'etale topos by a natural map
\[ 
	\varepsilon \colon X_\et \to X_{{\rm n.}\et}.
\]
Theorem~\ref{thm:affine-kpi1} implies that for $\FF_p$-schemes we can compute \'etale cohomology of locally constant constructible sheaves using the naive \'etale topology.

\begin{corollary}
	Let $X$ be an $\FF_p$-scheme and let $\Ff$ be locally constant constructible sheaf on $X_\et$. Then the maps
	\[ 
		\varepsilon^*\colon H^*(X_{{\rm n.}\et}, \varepsilon_*\Ff)\to H^*(X_\et, \Ff)
	\]
	are isomorphisms. 
\end{corollary}

(Presumably the same assertion holds for general constructible sheaves under suitable finiteness conditions on $X$.)

Suppose that $X$ is a normal noetherian $\FF_p$-scheme. Then the naive \'etale site of $X$ is locally connected, and applying the Verdier functor one can associate to it the `naive \'etale' homotopy type $\Pi(X_{{\rm n.}\et})$. It follows that the natural map
\[ 
	\Pi(\varepsilon) \colon \Pi(X_{{\rm n.}\et}) \to \Pi(X_\et)
\]
is a $\natural$-isomorphism, and hence a weak equivalence.

By the results of Section~\ref{s:local-rigid}, analogous results hold for rigid spaces in positive and mixed characteristic, in which case the naive \'etale topology is generated by \emph{admissible} open coverings and finite \'etale covers of affinoids.

\subsection{Implications in \'etale homotopy}

Theorem~\ref{thm:an-kpi1} yields a ``finite'' description of the homotopy type of a smooth $n$-dimensional variety in characteristic $p$ in terms of the single profinite group $\pi_1(\AN n)$.  

\begin{corollary} \label{cor:homotopy}
	Let $(X, \bar x)$ be a pointed connected quasi-projective scheme over $k$, smooth of dimension $n$. Then there exist $n+1$ affine open subsets $U_0, \ldots, U_{n}\subseteq X$ containing $\bar x$ and covering $X$, and for every non-empty subset $I\subseteq \{0, \ldots, n\}$, a finite \'etale map 
	\[ 
		f_I \colon U_I \to \AN{n} \quad\text{where} \quad U_I = \bigcap_{i\in I} U_i
	\]
	with $f_I(\bar x)=0$. Thus each $U_I$ is a $K(\pi, 1)$, and $f_I$ induces an isomorphism of $\pi_1(U_I, \bar x)$ with an open subgroup $\Pi_I\subseteq \pi_1(\AN n, 0)$. For $I\subseteq J$, let $h_{IJ}\colon \Pi_J\to \Pi_I$ be the homomorphism induced by the inclusion $U_J\subseteq U_I$. Then the \'etale homotopy type of $X$ is the homotopy colimit of the diagram $(\{\mathbf{B}\Pi_I\}, \{h_{IJ}\})$.
\end{corollary}

In principle, this tells us that a good understanding of the group $\pi_1(\AN n)$ would shed light on the \'etale homotopy types of smooth $k$-schemes. Unfortunately, this group is too complicated for us to derive any concrete corollaries from the above presentation. 

In any case, the results seem to suggest that a very strong form of Grothendieck's anabelian conjectures could be true in positive characteristic. We allow ourselves to put forth some ambitious-looking questions in this direction in \S\ref{ss:anabelian-geometry}.

\subsection{Examples and complements}

We finish the paper with a few examples: 
\begin{itemize}
	\item Example~\ref{ex:linear}, showing that in the presence of fierce ramification, linear projections are not enough in general in the context of Theorem~\ref{thm:bertini}.
	\item Example~\ref{ex:hyperplane} of a smooth affine variety $X$ (the complement of a hyperplane arrangement) over $\ZZ$ such that $X_{\bar\FF_p}$ is a $K(\pi, 1)$ for every $p>0$ while $X_{\bar\QQ}$ is not. 
	\item Example~\ref{ex:differentan} showing that $\pi_1(\AN{n})$ and $\pi_1(\AN m)$ are not isomorphic as pro-finite groups for $n\neq m$, even though they have the same finite quotients.
\end{itemize}

We also study some abstract properties of fundamental groups of affine schemes in \S\ref{ss:pro-p}, comment on the relationship between our work and the `$K(\pi, 1)$ pro-$\ell$' neighborhoods of Friedlander and Gabber in \S\ref{ss:kpi1-pro-ell}, and state some open questions in the spirit of anabelian geometry in \S\ref{ss:anabelian-geometry}.

\subsection{Outline}

In Sections \ref{s:wild}--\ref{s:bertini} and Appendix~\ref{s:appendix}, we deal with the proof of the Bertini theorem. We start with a review of relevant ramification theory in Section~\ref{s:wild}. Then we prove the easy case of the theorem when the sheaf is non-fiercely ramified at infinity in \S\ref{ss:non-fierce}, and proceed to the general case in \S\ref{ss:bertini}. Appendix~\ref{s:appendix} contains an alternative proof of the rank one case of Theorem~\ref{thm:bertini}.

In Sections \ref{s:kpi1}--\ref{s:local-rigid}, we deal with $K(\pi, 1)$ schemes and rigid analytic spaces. In Section~\ref{s:kpi1}, we review the notion of a $K(\pi, 1)$ scheme. Then in Section~\ref{s:affine-kpi1}, we prove Theorem~\ref{thm:affine-kpi1}. The subsequent Section~\ref{s:local-rigid} treats the mixed characteristic and rigid geometry analogues of Theorem~\ref{thm:affine-kpi1}.  

In the last Section~\ref{s:complements}, we provide relevant examples and further discussion as listed above.

\subsection{Acknowledgements}

I would like to thank Ahmed Abbes, Bhargav Bhatt, Ofer Gabber, Kiran Kedlaya, Laurent Lafforgue, Martin Olsson, Arthur Ogus, Fabrice Orgogozo, Takeshi Saito, Vasudevan Srinivas, and Karol Szumilo for helpful conversations. I am especially grateful to Takeshi Saito for his help with the proof of Theorem~\ref{thm:bertini}, and to Ofer Gabber for pointing out that Theorem~\ref{thm:affine-kpi1} follows from its special case Corollary~\ref{cor:small-kpi1}. We thank an anonymous referee for pointing out a mistake in an earlier version of the paper and for many valuable comments. We would also like to thank Maciej Borodzik for providing Figure~\ref{figure:varpi}. The author was supported by NCN OPUS grant number UMO-2015/17/B/ST1/02634.

\setcounter{tocdepth}{1}
\tableofcontents

\section{Review of wild ramification}
\label{s:wild}

Let $k$ be an algebraically closed field of characteristic $p>0$.

\subsection{The Swan conductor} 

Let $C$ be a smooth curve over $k$, $x\in C(k)$ a point, $\Ff$ a locally constant constructible $\FF_\ell$-sheaf on $C\setminus\{x\}$. The Swan conductor $\Swan_x(\Ff)$ is an integer measuring the wild ramification of $\Ff$ at $x$. It depends only on the restriction of $\Ff$ to $\Spec K$ where $K$ is the fraction field of the henselization of $\Oo_{C, x}$. It appears in the Grothendieck--Ogg--Shafarevich formula \cite[Exp.~X, formula~7.2]{SGA5}
\begin{equation} \label{eqn:gos} 
	\chi_c(C, \Ff) = \rk(\Ff)\cdot\chi_c(C, \FF_\ell) - \sum_{x\in S} \Swan_{x}(\Ff). 
\end{equation}
Here $C$ is a smooth geometrically integral curve with a smooth projective model $\bar C$, $S = \bar C\setminus C$, and $\Ff$ is a locally constant constructible $\FF_\ell$-sheaf on $C$. It is often more convenient to use the \emph{total dimension}, defined as
\[ 
	\dimtot_x (\Ff) = \rk (\Ff) + \Swan_x (\Ff)
\]
instead of the Swan conductor. For example, the version of the above formula without compact supports is
\begin{equation} \label{eqn:gos-no-support} 
	\chi(C, \Ff) = \rk(\Ff)\cdot\chi(C, \FF_\ell) - \sum_{x\in S} \dimtot_x \Ff. 
\end{equation}

Let $\eta = \Spec K$ where $K$ is a henselian discrete valuation field containing $k$ whose residue field $\kappa$ is finitely generated over $k$. Let $\ell\neq p$ be a prime and let $\Ff$ be a locally constant constructible $\FF_\ell$-sheaf on $\eta$. Choose a separable closure $K^{\rm sep}$ of $K$, and let $\bar\eta = \Spec K^{\rm sep}$. Then $\Ff$ corresponds to the continuous ${\rm Gal}(K^{\rm sep}/K)$-module $M=\Ff(\bar\eta)$.

\begin{definition} \label{def:swan}
	\begin{enumerate}[(1)]
		\item We call $\Ff$ \emph{non-fiercely ramified} if there exists a finite separable Galois extension $K'/K$ such that the pullback of $\Ff$ to $\Spec K'$ is constant and such that the residue field extension $\kappa'/\kappa$ is separable. 
		\item Let $K'/K$ be as in (1), and let $x \in \Oo_{K'}$ be a generator of $\Oo_{K'}$ as an $\Oo_K$-algebra \cite[Chap.~III, Prop.~12]{SerreCorpsLocaux}. The \emph{ramification groups} $G_i \subseteq G = {\rm Gal}(K'/K)$ (cf. \cite[Chapter IV, \S 1]{SerreCorpsLocaux}) are defined as
		\[ 
			G_i = \{ \sigma \in G \, : \, \nu_{K'}(\sigma(x)-x) \geq i+1\}.
		\]

		They are independent of the choice of $x$. The group $G_0$ is the inertia subgroup of $G$, and $G_1$ is called the wild inertia subgroup.
		\item Suppose that $\Ff$ is non-fiercely ramified. Let $K'/K$ be as in (1), let $G={\rm Gal}(K'/K)$ (so that ${\rm Gal}(K^{\rm sep}/K)$ acts on $M$ through its quotient $G$), and let $G_i\subseteq G$ ($i\geq 0$) be the ramification groups.  The \emph{Swan conductor} of $\Ff$ is defined as (cf. \cite[\S 19.3]{SerreLinearReps}, \cite[\S 1.1]{Laumon})
		\[ 
			\Swan(\Ff) = \sum_{i=1}^\infty \frac{1}{[G_0:G_i]} \dim_{\FF_\ell} (M/M^{G_i}).
		\]
		It is an integer, and is independent of the choice of $K'$. We also define the \emph{total dimension}
		\[ 
			\dimtot(\Ff) = \dim_{\FF_\ell} M + \Swan(\Ff).
		\]
	\end{enumerate}
\end{definition}

\subsection{The Deligne--Laumon theorem}

In the context of Theorem~\ref{thm:bertini}, the utility of the Swan conductor comes from the following result of Deligne and Laumon.

\begin{theorem}[{\cite[Th\'eor\`eme 2.1.1]{Laumon}}] 
	Let $S$ be a noetherian excellent scheme, and let $f\colon X\to S$ be a separated morphism, smooth of relative dimension $1$. Let $Y\subseteq X$ be a closed subscheme which is finite and flat over $S$. Let $U=X\setminus Y$, let $\ell$ be a prime invertible on $S$, and let $\Ff$ be a locally constant constructible $\FF_\ell$-sheaf on $U$ of constant rank $r$. Consider the function $\phi\colon S\to \ZZ$ defined as follows
	\[ 
		\phi(s) = \sum_{y\in Y_{\bar s}} \left(\Swan_y(\Ff|_{U_{\bar s}}) + r\right)  = \sum_{y\in Y_{\bar s}} \dimtot_y (\Ff|_{U_{\bar s}})
	\]
	(here $\bar s$ is any geometric point over $s$, and the value of the function does not depend on the choice of $\bar s$). Then 
	\begin{enumerate}[(i)]
		\item The function $\phi$ is a constructible and lower-semicontinuous,
		\item If $\phi$ is locally constant on $S$, then the triple $(X, \Ff, f)$ is universally locally acyclic.
	\end{enumerate}
\end{theorem}

\begin{corollary}\label{cor:deligne-laumon}
	Let $f\colon X\to S$ be a projective morphism with geometrically connected fibers, smooth of relative dimension $1$, $i\colon S\to X$ a section, $\Ff$ a locally constant constructible $\FF_\ell$-sheaf on $U=X\setminus i(S)$. Suppose that the number $\Swan_{i(\bar s)}(\Ff|_{U_{\bar s}})$ is independent of the geometric point $\bar s$ of $S$. Then the sheaves $R^q f_* \Ff$ and $R^q f_! \Ff$ are locally constant with formation commuting with base change for all $q\geq 0$. In particular, we have $R^q f_* \Ff = 0$ for $q>1$.
\end{corollary}

\begin{proof}
We note first that for a locally constant constructible $\FF_\ell$-sheaf on $\eta=\Spec K$ where $K$ is a henselian discrete valuation field with perfect residue field, we have $\Swan(\Ff) = \Swan(\Ff^\vee)$. Indeed, it is clear from the fact that (using the notation of Definition~\ref{def:swan}) the $G_i$ are $p$-groups for $i\geq 1$, and hence $M$ is semisimple as a $G_i$-representation (by Maschke's theorem), so $\dim_{\FF_\ell} (M/M^{G_i}) = \dim_{\FF_\ell} (M^\vee/(M^\vee)^{G_i})$ and $\Swan(M) = \Swan(M^\vee)$.

It follows that we can apply \cite[Corollaire 2.1.2]{Laumon} (together with \cite[Remarque 2.1.3]{Laumon}) to both $\Ff^\vee$ and $\Ff$ to see that $R^q f_! \Ff^\vee$ and $R^q f_! \Ff$ are locally constant for $q\geq 0$. By Poincar\'e--Verdier duality, the sheaves $R^q f_* \Ff$ are then locally constant with formation commuting with base change.
\end{proof}

\subsection{The characteristic cycle of a constructible sheaf}

The recent work of Beilinson \cite{Beilinson} and Saito \cite{Saito2016} provides an analogue of the classical theory of the singular support and the characteristic cycle \cite[Chapter IX]{KashiwaraShapira} for constructible \'etale sheaves, fulfilling an expectation of Deligne. Let us review the relevant points briefly, following \cite{Saito2016}.

Let $X$ be a smooth scheme over $k$ which is everywhere of dimension $n$, let $\ell\neq p$ be a prime, and let $\Ff$ be a constructible complex of $\FF_\ell$-vector spaces on $X$. In \cite{Beilinson}, Beilinson defines the singular support $SS\Ff$ inside the cotangent bundle $T^* X$. It is the smallest closed conical subset $C\subseteq T^* X$ such that $\Ff$ is micro-supported on $C$ (cf. Definition~\ref{def:c-transversal}(4) below). He proves that all of its irreducible components have dimension $n$.

T.~Saito \cite{Saito2016} employed Deligne's ideas to define the characteristic cycle $CC\Ff$. It is an integral combination of the irreducible components of $SS\Ff$. It is uniquely determined by the property of being compatible with \'etale base change and by the Milnor formula 
\[ 
	-\dimtot \phi_x (\Ff, f) = (CC\Ff, df)_{T^* X, x}
	\quad
	\text{(cf. \cite[5.15]{Saito2016})}
\]
for the total dimension of the vanishing cycles of a morphism $f\colon X\to Y$ is to a smooth curve which is $SS\Ff$-transversal away from $x$ (cf. Definition~\ref{def:c-transversal}(2)).

If $X$ is a curve, $j\colon U\subseteq X$ a dense open subset, and $\Ff$ is a locally constant constructible $\FF_\ell$-sheaf on $U$, then
\[ 
	CC(j_! \Ff) = -\left(\rk \Ff \cdot [T^*_X X] + \sum_{x\in X\setminus U} \dimtot_x \Ff \cdot [T^*_x X]\right) 
\]
(cf. \cite[Lemma~5.11.3]{Saito2016}). Here $T^*_X X$ denotes the zero section in $T^* X$ and $T^*_x X$ the fiber at $x$.

\begin{definition} \label{def:c-transversal} 
	Let $C\subseteq T^* X$ be a conical (i.e., stable under the $\mathbf{G}_m$-action) closed subset. 
	\begin{enumerate}
		\item A morphism $h\colon W\to X$ from a smooth $k$-scheme $W$ is called \emph{$C$-transversal} if for every $w\in W(k)$, 
		\[ 
			C \cap \ker\left(h^*\colon  T^*_{h(w)} X\to T^*_w W\right) = \{0\} \text{ or }\varnothing.
		\]
		In this case, we define 
		\[ 
			h^\circ C = {\rm im}\left( W\times_X C \to W\times_X T^* X \xto{dh} T^* W\right).
		\]
		\item A morphism $f\colon X\to Y$ to a smooth $k$-scheme $Y$ is called \emph{$C$-transversal} if for every $x\in X(k)$,
		\[ 
			\text{the preimage of }C\text{ under }f^*\colon T^*_{f(x)} Y \to T^*_x X = \{0\} \text{ or }\varnothing.
		\]
		\item A pair of morphisms $h\colon W\to X$, $f\colon W\to Y$ of smooth $k$-schemes is called \emph{$C$-transversal} if $h$ is $C$-transversal and $f$ is $h^\circ C$-transversal. 
		\item We say that a constructible complex of $\FF_\ell$-sheaves $\Ff$ on $X$ is \emph{micro-supported on $C$} if for every $C$-transversal pair of morphisms $h\colon W\to X$, $f\colon W\to Y$, $f$ is locally acyclic with respect to $h^* \Ff$ (cf. \cite[Th.Finitude, Definition~2.12]{ThFinitude}).
	\end{enumerate}
\end{definition}

In our proof of the Bertini theorem, we have to control the wild ramification of the restrictions to curves of a given sheaf $\Ff$. To this end, we need some compatibility of the characteristic cycle with pull-back. 

\begin{definition} \label{def:properly-c-transversal}
	In the situation of Definition~\ref{def:c-transversal}, suppose that every irreducible component of $C$ has dimension $n$. A morphism $h\colon W\to X$ from a smooth $k$-scheme $W$ which is everywhere of dimension $m$ is called \emph{properly $C$-transversal} if it is $C$-transversal and if every irreducible component of $W\times_X C$ has dimension $m$. In this situation, let $A=\sum m_a [C_a]$ be an integral combination of the irreducible components of $C$. We define
	\[ 
		h^! A = (-1)^{n-m} (dh)_* h^* A.
	\]
	Here $h^*$ denotes pull-back along $W\times_X T^* X\to T^* X$, and $(dh)_*$ means the push-forward along $dh\colon W\times_X T^*X \to T^* W$ in the sense of intersection theory.
\end{definition}

\begin{theorem}[{\cite[Theorem~7.6]{Saito2016}}] \label{thm:ccpullback}
	Suppose that $h\colon W\to X$ is properly $SS\Ff$-transversal. Then
	\[ 
		CCh^* \Ff = h^! CC\Ff.
	\]
\end{theorem}

\begin{corollary} \label{cor:t-saito}
	Let $X$ be a smooth $k$-scheme, $D\subseteq X$ a divisor, $\Ff$ a locally constant constructible $\FF_\ell$-sheaf on $U=X\setminus D$. Then there exists a dense open subset $T^\circ\in D\times_X \PP(T X)$ with the following property: if $x\in D$ and $L\subseteq T_x X$ is a line such that the corresponding point $(x, L)$ lies in $T^\circ$, and if $C, C'\subseteq X$ are smooth locally closed curves with $C\cap D = \{x\} = C'\cap D$ and $T_x C = L = T_x C'$, then
	\[ 
		\Swan_x (\Ff|_{C\setminus \{x\}}) = \Swan_x (\Ff|_{C'\setminus \{x\}}).
	\]
\end{corollary}

Here we follow the convention that $\PP(T X)$ parametrizes lines in the tangent bundle $TX$ (this is consistent with \cite{Saito2016}). Thus points $D\times_X \PP(T X)$ are identified with pairs $(x, L)$ of a point $x\in D$ and a tangent direction $L\subseteq T_x X$.

\begin{proof}
We can assume that $X$ is everywhere of dimension $n$, and that $\Ff$ is of constant rank. Let $j\colon U\to X$ be the inclusion. Since every irreducible component of $SS(j_! \Ff)$ has dimension $n$, while $D$ has dimension $n-1$, there exists a dense open subset $D^\circ\subseteq D$ such that $\dim (SS(j_! \Ff) \cap T^*_x X)\leq 1$ for $x\in D(k)$. Replace $D$ with $D^\circ$ and set 
\[ 
	T^\circ =  \left\{ (x, L)\in D\times_X \PP(T X) \, :\, \omega|_L  \neq 0 \text{ for all nonzero }\omega\in SS(j_! \Ff)\times_X \{x\} \right\}.
\]
Then $T^\circ$ is a dense open subset of $D\times_X \PP(T X)$. 

Suppose that $C\subseteq X$ is a locally closed curve with $C\cap D=\{x\}$ and $(x, T_x C) \in T^\circ$, as in the statement. We check that the inclusion $i\colon C\to X$ is properly $SS(j_!\Ff)$-transversal. The condition that $(x, T_x C)\in T^\circ$ is equivalent to the fact that 
\[ 
	SS(j_! \Ff) \cap \ker\left( T^*_x X \to T^*_x C\right) = \{0\}.
\]
This means that $i$ is $SS(j_!\Ff)$-transversal. The additional condition on the dimensions of the irreducible components of $C\times_X SS(j_! \Ff)$ is satisfied automatically.

Theorem~\ref{thm:ccpullback} implies now that
\[ 
	CC(i^* j_!\Ff) = i^! CC(j_! \Ff),
\]
while
\[ 
	CC(i^* j_! \Ff) = - \rk\Ff \cdot [T^*_C C] + (\rk\Ff + \Swan_x(\Ff|_{C\setminus\{x\}}))\cdot [T^*_x C].
\]
By definition of $i^!$, the coefficient of $[T^*_x C]$ in $i^! CC(j_! \Ff)$ depends only on the map $T^*_x X \to T^*_x C$.
\end{proof}

\begin{remark}
In his slightly earlier paper \cite{Saito}, predating Beilinson's ideas, Saito defined the characteristic cycle of a locally constant constructible sheaf in a neighborhood of the generic point of the boundary divisor $D$ using a different method, and studied its behavior upon restrictions to curves. Our Corollary~\ref{cor:t-saito} can also be deduced from \cite[Corollary 3.9.2]{Saito}.
\end{remark}

\section{Proof of the Bertini theorem} 
\label{s:bertini}

In this section, $k$ remains to denote a fixed algebraically closed field of characteristic $p>0$. The assertions of Proposition~\ref{prop:bertini-nonfierce} and Theorem~\ref{thm:bertini-general} remain valid over any infinite characteristic $p$ field.

\subsection{The non-fierce case of the Bertini theorem}
\label{ss:non-fierce}

As a warm-up, we show that a variant of Theorem~\ref{thm:bertini} holds for sheaves with non-fierce ramification at infinity. In contrast with the general case, it is possible to choose the automorphism $\phi$ to be a general linear automorphism.

Let $X$ be an integral smooth scheme over $k$, $D\subseteq X$ an irreducible smooth divisor, $U=X\setminus D$ its complement. Let $\Ff$ be a locally constant constructible $\FF_\ell$-sheaf on $U$. Let $X_{(\eta_D)}$ denote the localization of $X$ at the generic point $\eta_D$ of $D$ for the \'etale topology. Then $X_{(\eta_D)}\times_X U$ is the spectrum of the henselization of the fraction field of $X$ with respect to the discrete valuation given by $D$.

\begin{definition} \label{def:nonfierce-sheaf}
	We call $\Ff$ \emph{non-fiercely ramified along $D$} if the restriction of $\Ff$ to $X_{(\eta_D)}\times_X U$ is non-fiercely ramified in the sense of Definition~\ref{def:swan}(1), and if this is the case we write $\Swan_D(\Ff) = \Swan(\Ff|_{X_{(\eta_D)}\times_X U})$. 
\end{definition}

\begin{proposition}[{cf. \cite[\S 2.2]{LaumonSurface}}] \label{prop:nonfierce-is-nice}
	In the above situation, suppose that $\Ff$ is non-fiercely ramified along $D$. Then there exists a dense open $D^\circ \subseteq D$ with the property that for any $x\in D^\circ(k)$ and any smooth locally closed curve $C\subseteq X$ with $C\cap D = \{x\}$ and transverse to $D$ at $x$, we have $\Swan_x(\Ff|_{C\setminus \{x\}}) = \Swan_D(\Ff)$.
\end{proposition}

\begin{proof}
We include a direct proof (surely standard) because we were unable to find one in the literature (but see Remark~\ref{rmk:nonfierce-through-char-cycle} below). Let $\bar\eta_D$ be a geometric point above $\eta_D$. Setting $\eta = X_{(\bar\eta_D)}\times_X U$ puts us in the henselian situation described in \S\ref{s:wild}. Let $K'/K$ be as in Definition~\ref{def:swan}. Since the residue field $\kappa$ of $K$ is separably closed, while $\kappa'/\kappa$ is separable because of the non-fierceness assumption, we have $\kappa'=\kappa$. Thus $K'$ is a totally ramified extension of $K$. By \cite[I \S 6]{SerreCorpsLocaux}, there exists an Eisenstein polynomial $P\in \Oo_K[T]$ such that $\Oo_{K'} \isom \Oo_K[T]/(P)$, and the image of the variable $T$ is a uniformizer of $\Oo_{K'}$ under this isomorphism.  

Both the non-fierceness assumption on $\Ff$ and the assertion of the proposition are \'etale local in a neighborhood of $\bar\eta_D$. Spreading out the data described in the first paragraph to an \'etale neighborhood, we can assume that
\begin{enumerate}[(1)]
	\item $X=\Spec A$ is affine,
	\item $D$ is principal, its ideal generated by an element $\pi\in A$,
	\item there exists a polynomial $P\in A[T]$ whose image in $\hat A[T]$ ($\hat A = \varprojlim A/\pi^{n+1}$) is an Eisenstein polynomial, such that, setting $B=A[T]/(P)$, $Y=\Spec B$, then $Y$ is normal and finite over $X$, and $V=\Spec B[1/\pi]$ is an \'etale torsor under a finite group $G$ over $U = \Spec A[1/\pi]$, 
	\item the pull-back of $\Ff$ to $V$ is constant,
	\item the polynomial $P$ has the form
		\[
			P = T^r + a_1 T^{r-1} + \ldots + a_r, 
			\quad
			a_i= u_i\cdot \pi^{m_i}, \, u_i\in A^\times, \, m_i \geq 1,
		\] 
	\item for all $\sigma\in G$, there exists an integer $m(\sigma)$ and a unit $u(\sigma)\in B^\times$ such that 
		\[ 
			\sigma(T) - T = u(\sigma)\cdot T^{m(\sigma)}.
		\]
\end{enumerate}

Let $C\subseteq X$ be a smooth locally closed curve through a point $x\in D$ and transverse to $D$ at that point. Transversality means that $\pi$ maps to a uniformizer of $\Oo_{C, x}$.  Let $C' = Y\times_X C = \Spec \Oo_C[T]/(P)$. Since the image of $P$ in $\Oo_{C, x}[T]$ is Eisenstein (by condition (5) above), $T$ gives a uniformizer of $C'$ at the unique point $x'$ over $x$. Moreover, condition (6) implies that
\[ 
	\nu_{C'}(\sigma(T)-T) = \nu_D(\sigma(T)-T), 
	\quad
	\
	(\sigma \in G)
\]
We deduce that the $G$-Galois extensions $\Oo_{C',x'}/\Oo_{C,x}$ and $B/A$ induce the same ramification filtration on $G$. This implies the required assertion. 
\end{proof}

\begin{remark} \label{rmk:nonfierce-through-char-cycle}
	The following argument employing the characteristic cycle was pointed out to us by an anonymous referee. Observe that by Theorem~\ref{thm:ccpullback} the assertion of Proposition~\ref{prop:nonfierce-is-nice} holds under the weaker assumption that
	\begin{equation} \label{eqn:nice-ss}
		SS(j_! \Ff) \subseteq T^*_D X \cup T^*_X X
	\end{equation}
	holds generically along $D$. If $\pi:Y\to X$ is as in the above proof the normalization of $X$ in a finite \'etale Galois cover $V\to U$ trivializing $\Ff$, then
	\begin{equation} \label{eqn:ss-incl1}
		SS(j_! \Ff) \subseteq SS(\pi_* \underline{\FF}_{\ell,Y}).  
	\end{equation}
	On the other hand, \cite[Lemma~4.2.6]{Saito2016} implies that 
	\begin{equation} \label{eqn:ss-incl2}
		SS(\pi_* \underline{\FF}_{\ell,Y}) \subseteq \pi_\circ SS(\underline{\FF}_{\ell,Y}) = \pi_\circ T^*_Y Y
	\end{equation}
	(cf. \cite[Definition~3.7]{Saito2016} for the definition of $\pi_\circ$). Now if $\Ff$ is non-fiercely ramified along $D$, there exists $V\to U$ as above such that $\pi: E:=\pi^{-1}(D)_{\rm red}\to D$ is separable. Passing to a neighborhood of $D$, we can assume that $\pi:E\to D$ is \'etale, in which case the commutative diagram
	\[ 
		\xymatrix{
			0\ar[r] & (\pi|_{E})^* N_{D/X} \ar[d] \ar[r] & (\pi|_{E})^* (\Omega^1_{X/k}|_D) \ar[d]^{\pi^*} \ar[r] & (\pi|_{E})^* \Omega^1_{D/k} \ar[d]^\isom \ar[r] & 0 \\
			0\ar[r] & N_{E/Y} \ar[r] & \Omega^1_{Y/k}|_{E} \ar[r] & \Omega^1_{E/k} \ar[r] & 0 
		}
	\]
	shows that for $x\in D$,
	\[ \pi_\circ (T^*_Y Y)\times_X \{x\} = \bigcup_{y\colon \pi(y)=x} \ker(\pi^*: T^*_x X \to T^*_y Y) \subseteq T^*_D X \times_X \{x\}, 
	\]
	and hence 
	\begin{equation} \label{eqn:ss-incl3}
		\pi_\circ(T^*_Y Y) \subseteq T^*_D X\cup T^*_X X. 
	\end{equation}
	Combining \eqref{eqn:ss-incl1}--\eqref{eqn:ss-incl3} we get \eqref{eqn:nice-ss}.

	Note that in the classical complex analytic setting \eqref{eqn:nice-ss} is always satisfied because the irreducible components of the singular support of a holonomic $\mathscr{D}$-module are Lagrangian subvarieties of $T^* X$, and hence are the closures in $T^* X$ of conormal bundles $T^*_Z X$ of smooth locally closed subschemes $Z\subseteq X$.  
\end{remark}

\begin{proposition} \label{prop:bertini-nonfierce}
	The assertion of Theorem~\ref{thm:bertini} holds for a general linear automorphism $\phi$ of $\AN{n+1}$ if the sheaf $\Ff$ is non-fiercely ramified along the hyperplane at infinity. 
\end{proposition}

\begin{proof}
By Proposition~\ref{prop:nonfierce-is-nice}, there exists a dense open subset $H^\circ\subseteq H$ of the hyperplane at infinity with the property that for any line $L\subseteq \AN{n+1}$ which meets $H^\circ$ at infinity the Swan conductor $\Swan_\infty(\Ff|_L)$ is independent of $L$. Therefore if we take for $\pi\colon \AN{n+1}\to \AN{n}$ a linear projection along a line $L$ which meets $H^\circ$ at infinity, then function $y\mapsto \Swan_\infty(\Ff|_{\pi^{-1}(y)})$ will be constant. Hence the sheaves $R^q \pi_* \Ff$ will be locally constant with formation commuting with base change by the Deligne--Laumon theorem (Corollary~\ref{cor:deligne-laumon}).
\end{proof}

\begin{remark}
We expect that for a locally constant constructible $\FF_\ell$-sheaf $\Ff$ on $\AN{n+1}$, there exists an automorphism $\phi$ of the form
\[ 
	\phi(x_1, \ldots, x_{n+1}) 
	= (x_1 + x_{n+1}^{d_1}, \ldots, x_{n} + x_{n+1}^{d_n}, x_{n+1})
\]
such that $\phi^* \Ff$ is non-fierce at infinity. 
\end{remark}

\subsection{The general case of the Bertini theorem} 
\label{ss:bertini}

Before going into the proof, let us explain its main idea. In the non-fierce situation in the previous section, the Swan conductor $\Swan_x(\Ff|_{C\setminus \{x\}})$ of the restriction of a sheaf $\Ff$ to a curve meeting the boundary divisor $D$ transversally at a single point $x$ was independent of $C$ for $x$ in a dense open $D^\circ\subseteq D$. 

In the general case, the theory of the characteristic cycle (Corollary~\ref{cor:t-saito}) shows that the same assertion holds if the tangent space $T_x C \subseteq T_x X$ is a fixed element of a dense open subset $T^\circ$ of $D\times_X \PP(T X)$. This implies that for $x\in D$ and a tangent direction $L\subseteq T_x X$ such that $(x, L)\in T^\circ$, the number $\Swan_x(\Ff_{C\setminus \{x\}})$ is independent of $C$ as long as $T_x C = L$. 

It would therefore suffice to construct an $\AA^1$-fibration $\AN{n+1}\to \AN n$  whose fibers meet the hyperplane at infinity transversally with the same tangent direction. This is probably impossible, but we can produce such a fibration whose fibers are tangent to order two to the hyperplane at infinity and agree to sufficiently high order at that point. If $p>2$, taking the normalization of their preimages in a cyclic covering of degree two ramified along the hyperplane at infinity makes them transverse to the boundary, with the same tangent direction (this is another idea due to T.~Saito). This allows one to apply Corollary~\ref{cor:t-saito} to the cyclic covering.

\begin{theorem} \label{thm:bertini-general}
	Let $k$ be a field of characteristic $p$. Let $\ell \neq p$ be a prime and let $\Ff$ be a locally constant constructible sheaf of $\FF_\ell$-vector spaces on $\AN{n+1}$ ($n\geq 0$). Consider the map 
	\[ 
		\varpi\colon \AN{n+1}\to \AN{n}, \quad \varpi(x_1, \ldots, x_{n+1}) = (x_1 - x_{n+1}^2, x_2, \ldots, x_n).
	\]
	Let $G$ be the group of affine automorphisms $\psi$ of $\AN{n+1}$satisfying $\psi^* x_1 = x_1$. Then there exists a dense open $G^\circ(\Ff)\subseteq G$ such that for all $\psi\in G^\circ(\Ff)(\bar k)$, $\psi^* \Ff$ is well-aligned with respect to $\varpi$.
\end{theorem}

(We note that this implies Theorem~\ref{thm:bertini}: since $G$ has an open subset isomorphic to an open subset of the affine space, $G(k)$ is dense in $k$ as long as $k$ is infinite, and hence $G^\circ(\Ff)(k)\neq\emptyset$. Let 
\[
	\phi_0\colon \AN{n+1}\to\AN{n+1}, \quad \phi_0(x_1, \ldots, x_{n+1}) = (x_1 - x_{n+1}^2, x_2, \ldots, x_n, x_{n+1}).
\] 
Then $\varpi = \pi\circ \phi_0$, so if $\psi^* \Ff$ is well-aligned with respect to $\varpi$ then $\phi^* \Ff$ is well-aligned with respect to $\pi$ where $\phi = \psi\circ\phi_0^{-1}$.)

\begin{proof}
Let the coordinates on $\AN{n+1}$ be $x_1, \ldots, x_{n+1}$. We put $\AN{n+1}$ inside $\PP^{n+1}_k$ in the usual way, by adding an additional homogeneous coordinate $x_0$. We will look at what happens `at infinity' by considering the open subset 
\[ X = D^+(x_1) \subseteq \PP^{n+1}_k. \]
It is isomorphic to $\AN{n+1}$, with coordinates $z=x_0/x_1$, $z_i = x_i/x_1$ ($i=2, \ldots, n+1$). Let $D=V(z)\subseteq X$ be the intersection of $X$ with the hyperplane at infinity and let $U=X\setminus D$. 

Consider the cyclic cover of degree two of $X$ ramified along $D$:
\[ 
	\sigma \colon X' = \Spec \Oo_X[w]/(w^2 - z) \to X,
\]
and let 
\[ 
	U' = \sigma^{-1}(U), \quad D' = X'\setminus U' = \sigma^{-1}(D)_{\rm red}, \quad \Ff' = (\sigma|_{U'})^* \Ff. 
\] 
We apply Corollary~\ref{cor:t-saito} to the tuple $(X', D', U', \Ff')$, obtaining the subset $T^\circ(\Ff)\subseteq D'\times_{X'} \PP(T X')$.  We have $X'\isom \AN{n+1}$ with coordinates $w, z_2, \ldots, z_{n+1}$, and $\sigma(w, z_2, \ldots, z_{n+1}) = (w^2, z_2, \ldots, z_{n+1})$. The commutative diagram of exact sequences
\[ 
	\xymatrix{
		0\ar[r] & (\sigma|_{D'})^* N_{D/X} \ar[d]_0 \ar[r] & (\sigma|_{D'})^* (\Omega^1_{X/k}|_D) \ar[d]\ar[r] & (\sigma|_{D'})^* \Omega^1_{D/k} \ar[d]^\isom \ar[r] & 0 \\
		0\ar[r] & N_{D'/X'} \ar[r] & \Omega^1_{X'/k}|_{D'} \ar[r] & \Omega^1_{D'/k} \ar[r] & 0 
	}
\]
shows that there is a natural splitting
\[ 
	\Omega^1_{X'}|_{D'} \isom N_{D'/X'} \oplus \Omega^1_{D'/k}, 
\]
using which we can regard $\PP(T_{D'} X')$ as a subset of $D'\times_{X'} \PP(T X')$. We set
\[ 
	T^{\circ\circ} = \left(D'\times_{X'} \PP(T X')\right) \setminus \left(\PP(T_{D'} X')\cup \PP(T D')\right).
\]

Note that the group $G$ acts naturally on $\PP^{n+1}_k$ as the group of automorphisms fixing both $x_0$ and $x_1$. Consequently, the action of $G$ preserves not only $\AN{n+1}$ but also $X$, on which it acts as the group of affine transformations fixing the coordinate $z=x_0/x_1$. Furthermore, since as $z$ is fixed, the action of $G$ on $X$ lifts to an action on $X'$ fixing $D'$. We claim that $T^{\circ\circ}$ is the open orbit of the action of $G$ on $D'\times_{X'} \PP(T X')$. Writing $X=\AA^1_z \times V$, $X'=\AA^1_w\times V$ where $V$ is the affine space $\AA^n_{z_2, \ldots, z_{n+1}}$, the action of an element $g\in G$ can be presented as 
\[ g(z, v) = (z, Av + v_0 + zv_1), \quad g(w, v) = (w, Av + v_0 + w^2 v_1)   \]
for some $A\in GL(V)$, $v_0, v_1\in V$. Present an element $\xi$ of $D'\times_{X'} T X'$ in the form $\xi = (0, v) + \varepsilon\cdot (w, \tilde v)$ where $\varepsilon^2 = 0$, then 
\[ 
	g(\xi) = (0, Av+v_0) + \varepsilon\cdot(w, A\tilde v) + \varepsilon^2\cdot (0, w^2 v_1) = (0, Av+v_0) + \varepsilon\cdot(w, A\tilde v).
\]
On the other hand, the class of $\xi$ lies in $T^{\circ\circ}$ if and only if $w\neq 0$ and $\tilde v \neq 0$. We can thus assume $w=1$. For $\xi' = (0, v') + \varepsilon\cdot (1, \tilde v')$ whose class lies in $T^{\circ\circ}$, take $A\in GL(V)$ satisfying $A\tilde v = \tilde v'$ and set $v_0 = v' -Av$, $v_1 = 0$, then $g$ as above sends $\xi$ to $\xi'$.

We define $G^\circ(\Ff)\subseteq G$ as the subscheme of all $\psi\in G$ which map the point 
\[ 
	(x', L')\in D'\times_{X'} \PP(T X')\colon \quad x' = (0, \ldots, 0), \quad  L' = k\cdot(1,0,\ldots,0, 1)
\] 
to a point in $T^\circ(\Ff)$. (This particular choice will be explained shortly by the calculations of the fibers of $\varpi$.) Since $T^\circ(\Ff)$ is a dense open and $(x', L')$ lies in the open orbit $T^{\circ\circ}$, $G^\circ(\Ff)$ is a dense open subscheme of $G$. Thus $\psi\in G^\circ(\Ff)(k)$ means that the point $(x', L')$ defined above lies in $T^\circ(\psi^* \Ff)$. 

We analyze the fibers of $\varpi$ (cf. Figure~\ref{figure:varpi}). The fiber over $y=(y_1, \ldots , y_n)\in \AN{n}(\bar k)$ can be parametrized by
\[ 
	\gamma(t) = (y_1 + t^2, y_2, \ldots, y_n, t) \,:\, \AN{1}\to \AN{n+1}.
\] 
This extends to a map $\bar\gamma\colon \PP^1_k\to \PP^{n+1}_k$ as follows:
\[ 
	\bar\gamma(t:s) = (s^2: y_1 s^2 + t^2: y_2 s^2 : \ldots : y_n s^2 : ts).
\]
We want to see what happens at infinity, so we look at $\beta(s) = \bar\gamma(1:s)$:
\[ 
	\beta(s) = (s^2 : y_1 s^2 + 1 : y_2 s^2 : \ldots : y_n s^2 : s).
\]
So $\beta(0) = (0:1:0:\ldots:0) \in X$, which corresponds to the point $(z, z_2, \ldots, z_{n+1}) =(0, \ldots, 0)\in X$, which we denote henceforth by $x$. Thus in the $z$-coordinates, $\beta$ takes the form
\[ 
	\beta(s) = \left( \frac{s^2}{y_1 s^2 + 1}, \frac{y_2 s^2}{y_1 s^2 + 1} , \ldots , \frac{y_n s^2}{y_1 s^2 + 1} , \frac{s}{y_1 s^2 + 1}\right)
\] 
(defined for $s$ such that $y_1s^2 +1 \neq 0$). 

Let $C$ be the closure of the image of $\beta_y$ in $X$. Then $C\cap D = \{x\}$, $C$ meets $D$ with multiplicity two at $x$ and is smooth at $x$. Let $\tilde C=\AA^1\setminus\{ y_1 s^2 + 1 = 0 \} \to C$ be the normalization of $C$, $\tilde x=0$ the unique point of $\tilde C$ above $x$.

\emph{Case $p>2$:} The preimage of $C$ in $X'$ has two branches at $x'=(0,\ldots, 0)$, each one smooth and transverse to $D'$ (cf. top-right of Figure~\ref{figure:varpi}). This can be seen formally locally: look at
\[ 
		\beta'_{\pm}(s) = \left( \frac{\pm s}{\sqrt{y_1 s^2 + 1}}, \frac{y_2 s^2}{y_1 s^2 + 1} , \ldots , \frac{y_n s^2}{y_1 s^2 + 1} , \frac{s}{y_1 s^2 + 1}\right)
\]
where $\sqrt{y_1 s^2 + 1} = \sum \binom{1/2}{i} (y_1 s^2)^i$. So $\beta'_{\pm} \colon \Spf k[[s]] \to X'$ are two formal curve germs at $x'$, smooth and transverse to $D'$ at $x'$. Moreover, they map isomorphically to the normalization of the germ of $C$ at $x$ --- in other words, the diagram 
\[ 
	\xymatrix{
		\Spf k[[s]] \ar[d]_\cong \ar[rr]^{\beta'_\pm} & & \Spf \hat\Oo_{X', x'} \ar[d]^\sigma \\
		\Spf \hat\Oo_{\tilde C, 0}\ar[r]_\beta & \Spf \hat\Oo_{C, x} \ar@{^{(}->}[r] & \Spf \hat\Oo_{X, x}
	}
\] 
commutes. 

\emph{Case $p=2$:} In this case, since $(y_1 s^2 + 1) = (y_1^{1/2} s + 1)^2$, $\beta$ lifts (globally) along $\sigma$:
\[ 
 	\beta'_1(s) = \left( \frac{s}{y_1^{1/2} s + 1}, \frac{y_2 s^2}{y_1 s^2 + 1} , \ldots , \frac{y_n s^2}{y_1 s^2 + 1} , \frac{s}{y_1 s^2 + 1}\right), \quad \beta = \sigma\circ\beta'_1.
\]

In either case,  the tangent vector of $\beta'_1$ at $x'$ is
\[ 
	\frac{d}{ds} \beta'_{1}(s)|_{s=0} = (1, 0,\ldots, 0, 1) \in L'.
\]

\begin{figure}[ht]
	\centering
	\includegraphics[width=0.8\textwidth]{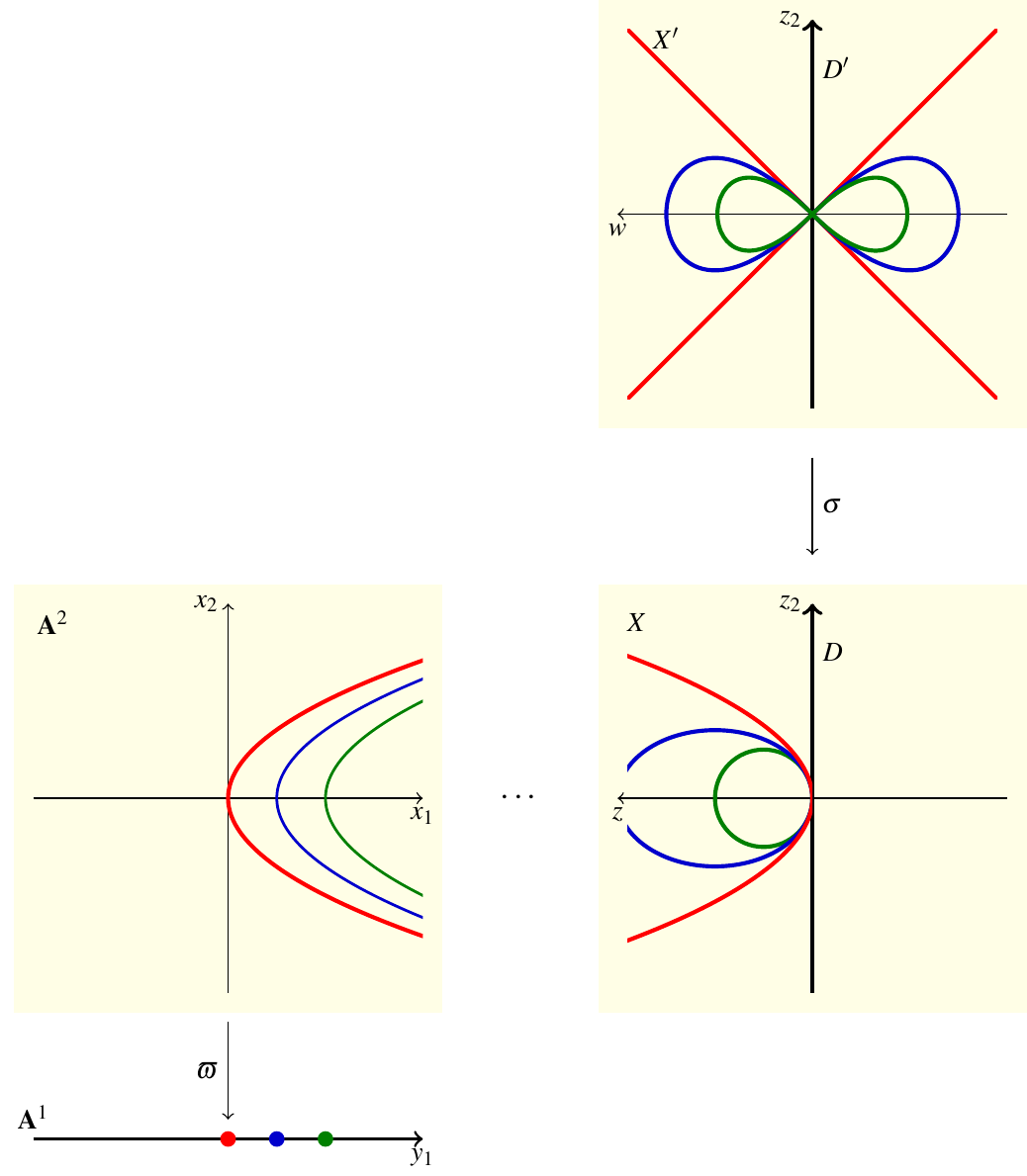} 
	\caption{Analysis of the fibers of $\varpi$ ($n=1$, $k=\mathbf{R}$).}
	\label{figure:varpi}
\end{figure}

Shrinking $X'$ to an \'etale neighborhood of $x'$, we can assume the existence of a locally closed curve $C'$ through $x'$ with germ $\beta'_1$ at $x'$ (more precisely, if $p>2$ , the \'etale cover comes from adjoining $\sqrt{y_1 s^2 + 1}$). Then $C'\to \tilde C$ is \'etale at $x'\mapsto \tilde x$, so by the fact that the Swan conductor depends only on the complete local field at $x$, we have 
\[ 
	\Swan_\infty(\Ff|_{\varpi^{-1}(y)}) = \Swan_x(\Ff|_{C\setminus\{x\}}) = \Swan_{x'}(\Ff|_{C'\setminus\{x'\}}).
\]

The tangent direction $(x', L')$ of $C'$ (which does not depend on $y$) lies in $T^\circ(\psi^* \Ff)$. The defining property of $T^\circ(\psi^* \Ff)$ implies thus that $\Swan_\infty(\Ff|_{\varpi^{-1}(y)})$ is independent of $y$. By the Deligne--Laumon theorem (Corollary~\ref{cor:deligne-laumon}), this implies that the sheaves $R^q \varpi_* (\psi^* \Ff)$ are locally constant, with formation commuting with base change.
\end{proof}

\section{\texorpdfstring{Review of $K(\pi, 1)$ schemes}{Review of K(pi, 1) schemes}} \label{s:kpi1}

\stepcounter{subsection}

Following Abbes and Gros \cite[\S 9]{Higgs2}, we will only consider schemes $X$ satisfying the following condition:
\begin{equation} \label{cond:qcqs}
	X\text{ \emph{is coherent and has finitely many connected components.}}
\end{equation}
By definition, a scheme is coherent if it is quasi-compact and quasi-separated. Note that if $A$ is a ring with finitely many idempotent elements, then $X=\Spec A$ satisfies \eqref{cond:qcqs}, and that if $X'\to X$ is a finite \'etale morphism and $X$ satisfies \eqref{cond:qcqs}, then so does $X'$.

\begin{definition}[{cf. \cite[Definition~9.20]{Higgs2}}] \label{def:kpi1}
	A pointed connected scheme $(X, \bar x)$ satisfying \eqref{cond:qcqs} is called a $K(\pi, 1)$ if for every locally constant constructible abelian sheaf $\Ff$ on $X$, the natural maps
	\begin{equation} \label{eqn:rhomaps2}
		\rho^q \colon H^q(\pi_1(X, \bar x), \Ff_{\bar x})\to H^q(X, \Ff)   
	\end{equation}
	are isomorphisms for all $q\geq 0$. This notion is independent of the choice of the geometric base point $\bar x$. We call a scheme $X$ satisfying \eqref{cond:qcqs} a $K(\pi, 1)$ if its connected components are $K(\pi, 1)$ schemes.
\end{definition}

Note that this is stronger than the notion used in \emph{op.cit.}, as we do not require that $\Ff$ be of torsion order invertible on $X$. The reference \cite[\S 9]{Higgs2} contains the most detailed discussion of this and related notions.

\begin{proposition} \label{kpi1prop} 
	Let $X$ be a scheme satisfying \eqref{cond:qcqs}. 
	\begin{enumerate}[(a)]
		\item $X$ is a $K(\pi, 1)$ if and only if for every locally constant constructible abelian sheaf $\Ff$ on $X$, and every class $\zeta\in H^q(X, \Ff)$ with $q>0$, there exists a finite \'etale surjective map $f\colon X'\to X$ such that $f^*(\zeta) = 0 \in H^q(X', f^* \Ff)$.
		\item Let $f\colon X'\to X$ be a finite \'etale surjective map. Then $X$ is a $K(\pi, 1)$ if and only if $X'$ is.
		\item $X$ is a $K(\pi, 1)$ if and only if for every prime $\ell$, every locally constant constructible $\FF_\ell$-sheaf $\Ff$ on $X$, and every class $\zeta\in H^q(X, \Ff)$ with $q>0$, there exists a finite \'etale surjective map $f\colon X'\to X$ such that $f^*(\zeta) = 0 \in H^q(X', f^* \Ff)$.
	\end{enumerate}
\end{proposition}

\begin{proof}
Assertions (ab) are \cite[Proposition 3.2(ab)]{Achinger}. For (c), note that functoriality of the long cohomology exact sequence implies that if
\[ 
	0\to \Ff' \to \Ff\to \Ff''\to 0
\]
is a short exact sequence of locally constant constructible sheaves and $\Ff'$ and $\Ff''$ satisfy the assertion of (a) after pulling back along every finite \'etale $f\colon X'\to X$, then so does $\Ff$. Since every locally constant constructible sheaf $\Ff$ has a finite filtration whose quotients are $\FF_\ell$-sheaves for various primes $\ell$, and every $X'$ finite \'etale over $X$ is a $K(\pi, 1)$ if $X$ is (by (b)), the `if' part of (c) follows, and the `only if' part is obvious.
\end{proof}

\begin{lemma} \label{lemma:kpi1q1}
	The maps $\rho^q$ are isomorphisms for $q\leq 1$. (In fact, this holds for $q=0$ and sheaves of sets, and for $q=1$ and sheaves of groups as well.) Therefore schemes of cohomological dimension $\leq 1$ (in particular, affine schemes of finite type of dimension $\leq 1$ over a separably closed field) are $K(\pi, 1)$.
\end{lemma}

\begin{proof}
For $q=0$, cf. \cite[Proposition 9.17]{Higgs2}. The statement for $q=1$ follows from the torsor interpretation of $H^1$ (cf. \cite[Remark~5.2]{Olsson}): a class $\zeta\in H^1(X, \Ff)$ corresponds to an isomorphism class of an $\Ff$-torsor $f\colon X'\to X$. The pullback $X'\times_X X'\to X'$ has a section, and hence is a trivial $f^*\Ff$-torsor, thus the corresponding class $f^* \zeta \in H^1(X', f^*\Ff)$ is zero. 
\end{proof}

\begin{proposition} \label{prop:kpi1artinmazur}
	Let $(X, \bar x)$ be a pointed connected noetherian scheme. Assume moreover that $X$ is geometrically unibranch (\cite[6.15.1]{EGAIV}, e.g. $X$ normal). Then $X$ is a $K(\pi, 1)$ if and only if $\pi_q(X,x)=0$ for $q>1$, where $\pi_q(X,x)$ is the \'etale homotopy group of Artin--Mazur \cite{ArtinMazur}.
\end{proposition}

\begin{proof}
Consider the natural map of sites $\rho\colon X_\et \to B\pi_1(X, x)$ and the associated map $\Pi \rho$ where $\Pi$ is the Verdier functor \cite[\S 9]{ArtinMazur}. On the one hand, $X$ is a $K(\pi, 1)$ if and only if $\Pi\rho$ is a $\natural$-isomorphism (cf. \cite[Theorem 4.3]{ArtinMazur}). On the other hand, $\Pi\rho$ induces an isomorphism on $\pi_1$ and $\pi_q(B\pi_1(X, x))=0$ for $q>1$, so $\pi_q(X,x)=0$ for $q>1$ if and only if $\Pi\rho$ is a weak equivalence. Both source and target of $\Pi\rho$ being pro-finite (thanks to $X$ being geometrically unibranch, \cite[Theorem~11.1]{ArtinMazur}), we conclude by \cite[Corollary 4.4]{ArtinMazur}.
\end{proof}

\begin{example} \label{ex:padickpi1}
Let $p$ be a prime. Then every connected affine $\FF_p$-scheme $X$ is a $K(\pi, 1)$ for $p$-torsion coefficients, that is, condition (c) of Proposition~\ref{kpi1prop} holds for $\ell=p$. Let $\Ff$ be a locally constant constructible $\FF_\ell$-sheaf on $X$, and let $\zeta\in H^q(X, \Ff)$ ($q>0$). We need to find a finite \'etale cover $X'\to X$ killing $\zeta$. First, we can assume that $\Ff$ is constant, as there exists a finite \'etale cover $X'\to X$ such that the pullback of $\Ff$ to $X'$ is constant, with $X'$ affine and connected. Second, we can reduce to the case $\Ff = \FF_p$. In this case, the Artin--Schreier sequence on $X_\et$
\[  0\to \FF_p \to \Oo_X \xlongrightarrow{1-F} \Oo_X \to 0 \]
together with Serre vanishing ($H^q(X_\et, \Oo_X) = H^q(X, \Oo_X) = 0$ for $q>0$) shows that $H^q(X, \FF_p) = 0$ for $q>1$. Thus if $q>1$, we are done. If $q=1$, then $\zeta$ corresponds to an $\FF_p$-torsor on $X_\et$, which again can be made trivial by a finite \'etale $X'\to X$. 

This example has been recently used by Scholze \cite[Theorem~4.9]{ScholzePAdic} to show that any Noetherian affinoid adic space over $\Spa(\QQ_p, \ZZ_p)$ is a $K(\pi, 1)$ for $p$-adic coefficients. We will follow Scholze's argument to prove that such spaces are in fact $K(\pi, 1)$ for all coefficients in Section~\ref{s:local-rigid}.
\end{example}

\section{\texorpdfstring{Affine $\FF_p$-schemes are $K(\pi, 1)$}{Affine Fp-schemes are K(pi, 1)}} \label{s:affine-kpi1}

\subsection{The affine space is a $K(\pi, 1)$}
\label{ss:an-kpi1}

We start by showing that $\AN n$ over a field $k$ of characteristic $p>0$ is a $K(\pi, 1)$ by induction on $n$, using Theorem~\ref{thm:bertini} in the induction step. Let us sketch the idea of the proof. By the characterization of Proposition~\ref{kpi1prop}(b), being a $K(\pi, 1)$ means being able to kill nonzero degree cohomology classes of locally constant constructible $\FF_\ell$-sheaves ($\ell$ an arbitrary prime) using finite \'etale covers. The case $\ell=p$ follows easily from Artin--Schreier theory, so suppose $\ell\neq p$. If $\Ff$ is such a sheaf, then Theorem~\ref{thm:bertini} for $\Ff$ implies that for a certain fibration $\pi\colon \AN{n+1} \to \AN{n}$, the higher direct image sheaves $R^i \pi_* \Ff$ are locally constant, and hence one can kill their cohomology using finite \'etale covers of $\AN{n}$. We derive the corresponding statement for $\Ff$ using the Leray spectral sequence of $\pi$.

\begin{theorem} \label{thm:an-kpi1}
	Let $k$ be a field. Then the affine space $\AN n$ is a $K(\pi, 1)$ scheme.
\end{theorem}

\begin{proof}
We prove this by induction on $n\geq 0$. We can assume that $k$ is an infinite field of characteristic $p>0$. Let $\Ff$ be a locally constant constructible abelian sheaf on $\AN{n+1}$. We want to show that for every class $\zeta\in H^q(\AN{n+1}, \Ff)$ ($q>0$) there exists a finite \'etale surjective $f\colon X'\to \AN{n+1}$ such that $f^*(\zeta)=0\in H^q(X', f^* \Ff)$. This is automatic for $q=1$ (Lemma~\ref{lemma:kpi1q1}), so we can assume $q>1$. Moreover, by Proposition~\ref{kpi1prop}(c), we can assume that $\Ff$ is an $\FF_\ell$-sheaf for a certain prime $\ell$. The case $\ell=p$ is handled by Example~\ref{ex:padickpi1}, so we can assume $\ell\neq p$.

By Theorem~\ref{thm:bertini}, there exists an $\AA^1$-bundle $\pi\colon \AN{n+1}\to\AN{n}$ such that the sheaves $R^i \pi_* \Ff$ are locally constant for $i\geq 0$, with formation commuting with base change. Since the fibers are $\AN 1$, we have $R^i \pi_*\Ff = 0$ for $i>1$, so the Leray spectral sequence for $\pi$ has only two nonzero rows, yielding an exact sequence
\[    
	\ldots \to H^q(\AN{n}, \pi_* \Ff) \to H^q(\AN{n+1}, \Ff) \to H^{q-1}(\AN{n}, R^1 \pi_* \Ff) \to \ldots. 
\]
Let $\zeta_0$ be the image of $\zeta$ in $H^{q-1}(\AN{n}, R^1 \pi_* \Ff)$. Since $q>1$ and $R^1 \pi_* \Ff$ is locally constant, the $K(\pi, 1)$ property of $\AN{n}$ implies that there exists a finite \'etale surjective $X\to \AN{n}$ killing $\zeta_0$. 

Replace $\AN{n+1}$ with $X' = X\times_{\AN{n}} \AN{n+1}$, $\pi\colon \AN{n+1}\to \AN{n}$ with its pullback $\pi'\colon X' \to X$, $\Ff$ and $\zeta$ with their pullbacks $\Ff'$, $\zeta'$ to $X'$. We again have an exact sequence as above, but now since $\zeta'$ maps to 0 in $H^{q-1}(X, R^1 \pi'_* \Ff')$, it is the pullback of a class $\zeta_1 \in H^q(X, \pi'_* \Ff')$. Again, since $X$ is a $K(\pi, 1)$ and $\pi'_* \Ff'$ is locally constant, we conclude that there is a finite \'etale surjective $Y\to X$ killing $\zeta_1$, and then $Y'=Y\times_X X'\to \AN{n+1}$ kills $\zeta$, as desired.
\end{proof}

\subsection{\'Etale schemes over the affine space} 
\label{ss:nagata}

Next, we deal with affine schemes endowed with an \'etale map to $\AN n$. To this end, we employ the following result.

\begin{proposition} \label{prop:nagatatrick}
	Let $k$ be a field of characteristic $p>0$. Let $U$ be an affine scheme of finite type over $k$, and let $g\colon U\to \AN{n}$ be an \'etale map. Then there exists a finite \'etale map $f\colon U\to \AN{n}$.
\end{proposition}

\begin{proof} 
This is a variant of Nagata's proof of Noether normalization (cf. \cite[I \S 1]{Mumford}). Write $U=\Spec R$, $R=k[x_1, \ldots, x_n, x_{n+1}, \ldots, x_r]/I$, where $x_1, \ldots, x_n$ are the pull-backs of the coordinates on $\AN{n}$ via $g$. We shall prove a slightly stronger statement: \emph{given any $x_1, \ldots, x_r\in R$ such that $r\geq n = \dim R$ and $x_1, \ldots, x_n$ are algebraically independent over $k$, there exist $y_1, \ldots, y_n \in R$ such that $R$ is finite over $k[x_1 + y_1^p, \ldots, x_n + y_n^p]$.} This implies what we want to prove because $dy_i^p = 0$, so $f=(x_1+y^p_1, \ldots, x_n+y^p_n)$ is \'etale if and only if $g$ is. 

The proof of this assertion is by induction on $r-n$: if $r=n$, $g$ is a closed immersion, and we take $f=g$. For the induction step, pick positive integers $a_1, \ldots, a_{r-1}$ and consider the elements
\[ 
	x'_i = x_i - x_r^{pa_i}, \quad i = 1, \ldots, nr-1. 
\]
Pick a nonzero $g\in I$, so that we have the relation
\[ 
	g(x'_1 + x_r^{pa_1}, \ldots, x'_{r-1} + x_r^{pa_{r-1}}, x_r) = 0 \text{ in }R.   
\]
By the usual argument, for $0 \ll a_1 \ll a_2 \ll \ldots \ll a_{r-1}$, this will be a monic polynomial in $x_r$ with coefficients in $k[x'_1, \ldots, x'_{r-1}]$. This means that if $R'\subseteq R$ is the subring generated by $x'_1, \ldots, x'_{r-1}$, then $x_r$ is integral over $R$. Since $x_i = x'_i + x_r^{pa_i}$, the other $x_i$ are also integral over $R'$, hence $R$ is integral over $R'$. As $R'$ is generated by $r-1$ elements, with $x'_1, \ldots, x'_n$ among them, we can apply the induction assumption to $R'$, $x'_1, \ldots, x'_n$ to find $y'_i \in R'$ such that $R'$ is finite over $k[x'_1+(y'_1)^{p}, \ldots, x'_n + (y'_n)^{p}]$. Thus $R$ is finite over $k[x_1+y^p_1, \ldots, x_n+y^p_n]$ where $y_i = y'_i - x_r^{a_i}$, 
\end{proof}

\begin{remark}
A related result has been obtained by Kedlaya \cite{Kedlaya}. As far as the author can tell, the trick of adding $p$-th powers to make a given map finite while preserving \'etaleness goes back to Abhyankar \cite{Abhyankar}. The author learned this technique from Katz's lectures \cite{Katz}. We have previously used a variant of this fact in a mixed characteristic situation \cite[Proposition~5.4]{Achinger}. In \S\ref{ss:nagata-rigid}, we will give a rigid analytic variant of Proposition~\ref{prop:nagatatrick}.
\end{remark}

\begin{corollary} \label{cor:small-kpi1}
	Let $U$ be an affine scheme of finite type over $k$, and let $g\colon U\to \AN{n}$ be an \'etale map. Then $U$ is a $K(\pi, 1)$ and for every geometric point $\bar u$ of $U$, $\pi_1(U, \bar u)$ is isomorphic to an open subgroup of $\pi_1(\AN n, 0)$.
\end{corollary}

\begin{proof}
This follows directly from Proposition~\ref{prop:nagatatrick}, Theorem~\ref{thm:an-kpi1}, and Proposition~\ref{kpi1prop}(b). 
\end{proof}

\subsection{Henselian pairs and Gabber's theorem}
\label{ss:henselian-pairs}

Recall \cite[Chapitre XI]{RaynaudHenselian} that a \emph{henselian pair} is a pair $(A, I)$ consisting of a ring $A$ and an ideal $I\subseteq A$ such that for every \'etale $A$-algebra $A'$, the restriction map
\[ 
	\Hom_A(A', A)\to \Hom_A(A', A/I)
\]
is a bijection. If $(A, I)$ is any pair, the \emph{henselization} $(A, I)\to (A^h, I^h)$ of $(A, I)$ is the initial map of pairs to a henselian pair. Henselization exists, and can be constructed as the inductive limit of pairs $(B, J)$ indexed by \'etale $A$-algebras $B$ endowed with a map of $A$-algebras $\sigma\colon B\to A/I$, and $J=\ker(\sigma)$ (cf. \cite[\S 0]{Gabber}). Moreover $I^h = I\cdot A^h$.

Let $(A, I)$ be a henselian pair. We set $X=\Spec A$, $X_0=\Spec A/I$, $i\colon X_0\to X$ the closed immersion. By \cite[\S 1]{Gabber}, the restriction functors 
\[ 
	i^* \colon \left( \text{finite \'etale covers of }X\right) 
	\to
	\left( \text{finite \'etale covers of }X_0\right), 
\]
\[ 
	i^*\colon \left(\text{lcc sheaves on }X\right) 
	\to
	\left(\text{lcc sheaves on }X_0\right) 
\]
are equivalences. For higher cohomology, we have the following.

\begin{theorem}[{\cite[Theorem 1]{Gabber}}] \label{thm:gabber-affine-analog}
	In the above situation, for every torsion abelian sheaf $\Ff$ on $X$ the restriction maps
	\[ 
		i^* \colon H^q(X, \Ff) \to H^q(X_0, i^* \Ff)
	\]
	are isomorphisms for all $q\geq 0$.
\end{theorem}

\begin{corollary} \label{cor:henselian-pairs}
	In the above situation, $X$ is a $K(\pi, 1)$ if and only if $X_0$ is a $K(\pi, 1)$.
\end{corollary}

\begin{proof}
If $Y=\Spec B$ is a finite \'etale $X$-scheme and $J=IB$, then $(B, J)$ is a henselian pair, and hence the above statements hold for $(B, J)$. Let $\Ff$ be an lcc sheaf on $X$, $\zeta\in H^q(X, \Ff)\isom H^q(X_0, i^* \Ff)$ ($q>0$). To show that $X$ (resp. $X_0$) is a $K(\pi, 1)$ means that for every such pair $(\Ff, \zeta)$, there exists a finite \'etale surjective $Y\to X$ (resp. $Y_0\to X_0$) killing $\zeta$. By the above remarks, $X$ is thus a $K(\pi, 1)$ if and only if $X_0$ is.
\end{proof}

\subsection{The general case}
\label{ss:proof-affine-kpi1}

Finally, we deal with general connected affine $\FF_p$-schemes. 

\begin{lemma} \label{lemma:limit-kpi1}
	Suppose that $A_\infty$ is a ring with no nontrivial idemponents which is the union of a filtered family of subrings $A_i\subseteq A_\infty$ such that $X_i = \Spec A_i$ is a $K(\pi, 1)$. Then $X=\Spec A_\infty$ is a $K(\pi, 1)$.
\end{lemma}

\begin{proof}
Let $\Ff$ be a locally constant constructible sheaf on $X$, $\zeta\in H^q(X, \Ff)$ where $q>0$. We need to find a finite \'etale surjective $f\colon X'\to X$ such that $f^*(\zeta) = 0 \in H^q(X', f^*\Ff')$. There exists a sheaf $\Ff_i$ on $X_i$ for some $i$ such that $\Ff\isom \Ff_i\otimes_{A_i} A_\infty$. Since the limit is filtered, we can restrict ourselves to $A_j$ for $j\geq i$; let $\Ff_j = \Ff\otimes_{A_i} A_j$. The natural map
\[ 
	\varinjlim_{j\geq i} H^q(X_j, \Ff_j) \to H^q(X, \Ff)
\]
is an isomorphism, and hence there exists a $j$ such that $\zeta$ is the image of a class $\zeta_j\in H^q(X_j, \Ff_j)$. Since $X_j$ is a $K(\pi, 1)$, there exists a finite \'etale surjective $f_j \colon X'_j\to X_j$, but then the base change $f\colon X'=X'_j\times_{X_j} X\to X$ kills $\zeta$, as desired.
\end{proof}

\begin{theorem*}[\ref{thm:affine-kpi1}]
	Every connected affine scheme over $\FF_p$ is a $K(\pi, 1)$ scheme.
\end{theorem*}

\begin{proof}
Let $X=\Spec A$ be a connected affine scheme over $\FF_p$. For a finite subset $S\subseteq A$, let $\FF_p[S]$ be the subring of $A$ generated by $S$. Then $A$ is the union of all such rings $\FF_p[S]$. Therefore by Lemma~\ref{lemma:limit-kpi1} it suffices to treat the case when $A$ is generated over $\FF_p$ by a finite number of elements $a_1, \ldots, a_n$. These elements exhibit $X$ as a closed subscheme of $\AN n=\Spec P$, $P=k[x_1, \ldots, x_n]$ ($k=\FF_p$) defined by an ideal $I\subseteq P$. Let $X^h=\Spec P^h$ be the henselization of $\AN n$ along $X$. By definition, $P^h$ is the inductive limit of \'etale $P$-algebras $B$ endowed with a section over $X$, i.e. a $P$-algebra homomorphism $B/IB\to P/I=A$. By Corollary~\ref{cor:small-kpi1}, each $\Spec B$ is a $K(\pi, 1)$. By Lemma~\ref{lemma:limit-kpi1}, $X^h$ is a $K(\pi, 1)$. But $(P^h, IP^h)$ is a henselian pair and $P^h/IP^h = P/I = A$. Thus $X=\Spec A$ is a $K(\pi, 1)$ by Corollary~\ref{cor:henselian-pairs}.
\end{proof}

\section{Mixed characteristic and rigid analytic variants}
\label{s:local-rigid}

\subsection{Review of rigid geometry}

We recall the setup for rigid geometry in the sense of Raynaud, for which we follow \cite{AbbesRigid} (albeit we will only need the noetherian case, as opposed to the more general idyllic case treated in that book). An \emph{adic ring} is a complete and separated topological ring $A$ admitting a finitely generated ideal $I\subseteq A$ such that the ideals $I^n$ are open and form a basis of neighborhoods of $0\in A$. Such an $I$ is called an \emph{ideal of definition}. A formal scheme $\Xx$ is called \emph{adic} if it is locally of the form $\Spf A$ for an adic ring $A$. An \emph{admissible blow-up} is a morphism $\phi\colon \Xx'\to\Xx$ of finite type between adic formal schemes which is isomorphic to the blow-up of $\Xx$ of a finitely generated open ideal. The category $\RR$ of coherent rigid spaces is the localization of the category $\mathbf{S}$ of noetherian quasi-compact formal schemes and morphisms of finite type with respect to admissible blow-ups. If $\Xx$ is an object of $\mathbf{S}$, we denote the associated object of $\RR$ by $\Xx^\rig$. We call $\Xx^\rig$ the \emph{associated rigid space} of $\Xx$. A \emph{formal model} of a coherent rigid space $X$ is an object $\Xx$ of $\mathbf{S}$ together with an isomorphism $X\isom \Xx^\rig$. An \emph{affinoid rigid space} is a coherent rigid space admitting a formal model of the form $\Spf(A)$ for a noetherian adic ring $A$ which locally admits a principal ideal of definition. 

One has natural notions of finite and \'etale morphisms in $\RR$, which allow one to define the rig-\'etale topos $X_{\et}$ of a rigid space $X$. 

\subsection{The Gabber--Fujiwara theorem}
\label{ss:gabber-fujiwara}

Let $(A, I)$ be a noetherian henselian pair, $\hat A$ the $I$-adic completion of $A$. Our goal is to compare the cohomology of the affinoid rigid space $X=\Spf(\hat A)^\rig$ to the cohomology of the scheme $U=\Spec A \setminus V(I)$. The GAGA functor \cite[Chap. VII]{AbbesRigid} produces a morphism of topoi
\[ 
	\varepsilon\colon X_\et \to U_\et,
\]
inducing equivalences as in \S\ref{ss:henselian-pairs}
\[ 
	\varepsilon^* \colon (\text{finite \'etale covers of }U)
	\isomlong 
	(\text{finite \'etale covers of }X),
\]
\[ 
	\varepsilon^* \colon (\text{lcc sheaves on }U)
	\isomlong 
	(\text{lcc sheaves on }X).
\]

\begin{theorem}[{Gabber--Fujiwara, \cite[Corollary 6.6.3]{Fujiwara}}] \label{thm:gab-fuj}
	In the above situation, for every torsion abelian sheaf $\Ff$ on $U$ the maps
	\[ 
		\varepsilon^*\colon H^q(U, \Ff) \to H^q(X, \varepsilon^*\Ff)
	\]
	are isomorphisms for $q\geq 0$.
\end{theorem}

\subsection{Affinoid rigid spaces in characteristic $p$}

With the Gabber--Fujiwara theorem in place, we can easily deduce from Theorem~\ref{thm:affine-kpi1} its rigid analytic variant.

\begin{definition} \label{def:rigid-kpi1}
	We call a rigid space $X$ a \emph{$K(\pi, 1)$ space} if for every locally constant constructible \'etale sheaf $\Ff$ on $X$, and every class $\zeta \in H^q(X, \Ff)$ for $q>0$, there exists a finite \'etale surjective map $f\colon X'\to X$ such that $f^*(\zeta) = 0 \in H^q(X', f^* \Ff)$.
\end{definition}

\begin{lemma} \label{lemma:kpi1rigid}
	Let $(A, I)$ be a noetherian henselian pair, $\hat A$ the $I$-adic completion of $A$. Then $\Spf(\hat A)^\rig$ is a $K(\pi, 1)$ rigid space if and only if $\Spec(A)\setminus V(I)$ is a $K(\pi, 1)$ scheme.
\end{lemma}

\begin{proof}
This is analogous to the proof of Corollary~\ref{cor:henselian-pairs}. Let $f\colon V\to U$ be a finite \'etale cover, $B$ the normalization of $A$ in $V$, $J=IB$. Then $(B, J)$ is a noetherian henselian pair, and $V=\Spec B\setminus V(J)$. Moreover, if $Y=\Spf(\hat B)^\rig$, then $g\colon Y\to X$ is the finite \'etale cover corresponding to $f$ via the equivalence $\varepsilon^*$, and the square
\[ 
	\xymatrix{
		V_\et \ar[d]_{f} \ar[r]^{\varepsilon} & Y_\et \ar[d]^g \\
		U_\et \ar[r]_{\varepsilon} & X_\et
	}
\]
$2$-commutes. Consequently, if $\Ff$ is an abelian sheaf on $U_\et$, then the square
\[ 
	\xymatrix{
		H^*(V, f^* \Ff) \ar[r]^{\varepsilon^*} & H^*(Y, \varepsilon^* (f^* \Ff)) \\
		H^*(U, \Ff)\ar[r]_{\varepsilon^*} \ar[u]^{f^*} & H^*(X, \varepsilon^*\Ff) \ar[u]_{g^*}
	}
\]
commutes. By Theorem~\ref{thm:gab-fuj}, the horizontal arrows are isomorphisms if $\Ff$ is a torsion sheaf.

Let $\Ff$ be a locally constant constructible sheaf on $U$, $\zeta\in H^q(U, \Ff)\cong H^q(X, \varepsilon^* \Ff)$. If $U$ (resp. $X$) is a $K(\pi, 1)$ then there exists a finite \'etale surjective $f\colon V\to U$ (resp., $g\colon Y\to X$) such that $f^*(\zeta)=0$ (resp. $g^*(\zeta)=0$). The equivalence is thus clear in view of the preceding discussion.
\end{proof}

\begin{theorem} \label{thm:kpi1-affinoid-char-p}
	Let $X=\Spf(A)^\rig$ be an affinoid rigid space such that $pA=0$. Then $X$ is a $K(\pi, 1)$ space.
\end{theorem}

\begin{proof}
	 Let $I\subseteq A$ be an ideal of definition, and let $U=\Spec A \setminus V(I)$. Since $A$ locally admits a principal ideal of definition, the map $U\to \Spec A$ is affine, and hence $U$ is an affine scheme. Thus $U$ is a $K(\pi, 1)$ scheme by Theorem~\ref{thm:affine-kpi1}. Lemma~\ref{lemma:kpi1rigid} implies that $X$ is a $K(\pi, 1)$ rigid analytic space.
\end{proof}

\subsection{Affinoid rigid spaces in mixed characteristic} 

To deduce the mixed characteristic case, we use perfectoid spaces. Since they are by definition adic spaces, not rigid spaces, we have to consider $K(\pi, 1)$ adic spaces: we will simply call an adic space a $K(\pi, 1)$ if the condition of Definition~\ref{def:rigid-kpi1} holds. If $X$ is a rigid space, $X^{\rm ad}$ the associated adic space, then the \'etale topoi $X_\et$ and $X^{\rm ad}_\et$ are equivalent (cf. \cite[Proposition~2.1.4]{Huber}), and $X$ is a $K(\pi, 1)$ rigid space if and only if $X^{\rm ad}$ is a $K(\pi, 1)$ adic space. Our argument is analogous to \cite[Theorem~4.9]{ScholzePAdic}. 

\begin{proposition} \label{prop:kpi1-perfectoid}
	Let $(A, A^+)$ be a perfectoid algebra over a perfectoid field $K$. Then ${\rm Spa}(A, A^+)$ is a $K(\pi, 1)$ adic space.
\end{proposition}

\begin{proof}
	Suppose first that $K$ has characteristic $p$. Without loss of generality, we can assume that $K$ is the completed perfection of $\FF_p((t))$. In this case, $A^+$ is a $K^\circ$-algebra. By \cite[Lemma 6.13 (i)]{Scholze}, $(A, A^+)$ is the completion of a filtered direct limit of $p$-finite perfectoid affinoid $K$-algebras $(A_i, A^+_i=A^\circ_i)$. If we set $X_i ={\rm Spa}(A_i, A^+_i)$, then every finite \'etale cover of $X$ comes by base change from a finite \'etale cover of $X_i$. Thus if $\Ff$ is a locally constant constructible sheaf on $X$, we can assume that there exists an $i_0$ such that $\Ff$ is the base change of a locally constant constructible sheaf $\Ff_{i_0}$ on $X_{i_0}$. For $i\geq i_0$, let $\Ff_i$ be the pullback of $\Ff_{i_0}$ to $X_i$. The natural map
	\[ 
		\varinjlim_{i\geq i_0} H^q(X_i, \Ff_i) \to H^q(X, \Ff)
	\] 
	is an isomorphism (\cite[Corollary~7.18]{Scholze} combined with \cite[Corollary~7.19]{Scholze}). These remarks show that if all $X_i$ are $K(\pi, 1)$, so is $X$. Therefore it suffices to treat the case when $X$ is $p$-finite.

	By definition, being $p$-finite means being the completed perfection of an affinoid algebra $(B, B^+=B^\circ)$ topologically of finite type over $K$. For such algebras, ${\rm Spa}(B, B^+)$ is a $K(\pi, 1)$ by the noetherian case (Theorem~\ref{thm:kpi1-affinoid-char-p}). But ${\rm Spa}(A, A^+)\to {\rm Spa}(B, B^+)$ induces an equivalence of the \'etale topoi (\cite[Corollary~7.19]{Scholze}), and hence ${\rm Spa}(A, A^+)$ is a $K(\pi, 1)$.

	Finally, we treat the case when $K$ has characteristic zero. Let $X^\flat / K^\flat$ be the tilt. Then by \cite[Theorem~7.12]{Scholze}, the \'etale topoi of $X$ and $X^\flat$ are equivalent, and the same holds for all finite \'etale covers of $X$ and $X^\flat$ is a compatible way. Thus $X$ is a $K(\pi, 1)$ if and only if $X^\flat$ is. 
\end{proof}

\begin{theorem} \label{thm:kpi1affinoid-mixed}
	Let $X$ be (1) an affinoid noetherian adic space over ${\rm Spa}(\QQ_p, \ZZ_p)$ or (2) $X=\Spf(A)^\rig$ for a noetherian $p$-adic ring $A$. Then $X$ is a $K(\pi, 1)$ space.
\end{theorem}

\begin{proof}
By \cite[Proposition~2.1.4]{Huber}, (2) follows from (1). By the proof of \cite[Theorem~4.9]{ScholzePAdic}, there exists a system of finite \'etale covers $X_i$ of $X$ and an affinoid perfectoid $X_\infty$ over $\CC_p$ such that 
\[ 
	X_\infty\sim \varprojlim X_i \quad (\text{cf. \cite[Definition~7.14]{Scholze}}).
\]
Then $X_\infty$ is a $K(\pi, 1)$ by Proposition~\ref{prop:kpi1-perfectoid}, and we conclude that $X$ is a $K(\pi, 1)$ by \cite[Corollary~7.18]{Scholze}.
\end{proof}

\subsection{Application to $p$-adic Milnor fibers}

We can apply the Gabber--Fujiwara theorem once again, now in mixed characteristic, to go back to the henselian case.

\begin{theorem} \label{thm:henselian-mixed}
	Let $A$ be a noetherian $\ZZ_{(p)}$-algebra such that $(A, pA)$ is a henselian pair. Then $\Spec A[1/p]$ is a $K(\pi, 1)$ scheme.
\end{theorem}

\begin{proof}
By Lemma~\ref{lemma:kpi1rigid}, $\Spec A[1/p]$ is a $K(\pi, 1)$ if and only if the rigid space $\Spf(\hat A)^\rig$ is a $K(\pi, 1)$, where $\hat A$ is the $p$-adic completion of $A$. The latter is a $K(\pi, 1)$ by Theorem~\ref{thm:kpi1affinoid-mixed}.
\end{proof}

Fix a strictly henselian discrete valuation ring $V$ with residue field $k$ of characteristic $p>0$ and fraction field $K$ of characteristic zero. Let $s=\Spec k$, $\eta=\Spec K$, $\bar\eta=\Spec\bar K$ where $\bar K$ is an algebraic closure of $K$. For a scheme $X$ of finite type over $S=\Spec V$ and a geometric point $\bar x$ of $X_s$, the \emph{Milnor fiber} of $X$ at $\bar x$ is the scheme
\[ 
	M_{\bar x} = X_{(\bar x)}\times_S \bar\eta.
\]

\begin{corollary} \label{cor:milnor-fiber}
	Let $X$ be a $V$-scheme of finite type. Then for every geometric point $\bar x$ of $X_s$, the Milnor fiber $M_{\bar x}$ is a $K(\pi, 1)$ scheme. 
\end{corollary}

This allows us to strenghten the main result of \cite{Achinger}, removing the log smoothness hypothesis (but only in the case of $X^\circ = X$).

\begin{corollary} \label{cor:better-achinger}
	Let $X$ be a $V$-scheme of finite type. Consider the Faltings' topos $\tilde E$ of $X_{\bar\eta}\to X$ and the morphism of topoi
	\[ 
		\Psi \colon X_{\bar\eta,\et} \to \tilde E.
	\]
	Let $\Ff$ be a locally constant constructible abelian sheaf on $X_{\bar\eta,\et}$. Then $R^q\Psi_* \Ff = 0$ for $i>0$, and the natural maps
	\[ 
		H^*(\tilde E, \Psi_* \Ff)\to H^*(X_{\bar\eta,\et}, \Ff)
	\]
	are isomorphisms.
\end{corollary}

\begin{proof}
	Corollary~\ref{cor:milnor-fiber} shows that condition (B) of \cite[\S 1.2]{Achinger} holds. This implies condition (D) of \emph{op.cit.} which is exactly our assertion.
\end{proof}

\subsection{\'Etale affinoids over the polydisc}
\label{ss:nagata-rigid}

Interestingly, we have a rigid analytic variant of Proposition~\ref{prop:nagatatrick} (cf. \cite[Proposition~5.10]{Achinger}). 

\begin{proposition} \label{prop:nagata-rigid}
	Let $K$ be a complete discretely valued field whose residue field $k$ is of characteristic $p>0$, $U$ an affinoid rigid analytic space of finite type over $K$, $\BB^n$ the formal $n$-polydisc over $K$ (the rigid-analytic generic fiber of $\ANK{n}{\Oo_K}$),  $g\colon U\to \BB^n$ an \'etale morphism. Then there exists a finite \'etale morphism $f\colon U\to \BB^n$.
\end{proposition}

\begin{proof}
Let $V$ be the valuation ring of $K$. Let $R$ be the Tate algebra corresponding to $U$, let $R^\circ\subseteq R$ be the integral subring, and let $x_1, \ldots, x_n\in R^\circ$ be the pull-backs of the coordinates on $\BB^n$ via $g$. Pick $x_{n+1}, \ldots x_r\in R^\circ$ such that $x_1, \ldots, x_r$ form a set of topological generators of $R^\circ$ over $V$. This gives us a presentation
\[ 
	R^\circ \cong V\langle x_1, \ldots, x_n, x_{n+1}, \ldots, x_r\rangle / I, 
	\quad
	R \cong R^\circ\left[\frac{1}{\pi}\right]
\]
where $\pi\in V$ is a uniformizer. Let $\Omega_{g} = \hat\Omega_{R^\circ/V\langle x_1, \ldots, x_n\rangle}$. Since $g$ is \'etale, we have $\Hh^0_\rig(\Omega_g) = 0$, and hence $\Omega_{g}$ is killed by a power of $\pi$, say $\pi^{m}\cdot \Omega_{g} = 0$. By \cite[Proposition~5.4]{Achinger} applied to the reduction of $g$ modulo $\pi$ and $N=p^{m+1}$, there exist elements $\bar y_1, \ldots, \bar y_n \in \bar R$ which are polynomials in $p^{m+1}$-th powers of elements of $\bar R$, such that the map $\Spec \bar R\to \AN n$ given by $\bar x_1 + \bar y_1, \ldots, \bar x_n + \bar y_n$ is finite. Let $y_i$ be any lifts of $\bar y_i$ to $R^\circ$ which are polynomials in $p^{m+1}$-powers of elements of $R^\circ$. Arguing as in \cite[Lemma~5.9]{Achinger}, we see that the map $f\colon U\to \BB^n$ given by $x_1+y_1, \ldots, x_n+y_n$ is \'etale. It is also finite, because its reduction modulo $\pi$ is. 
\end{proof}

\section{Examples and complements} \label{s:complements}

\subsection{Linear projections do not suffice}
\label{ex:linear}

Let $k$ be an algebraically closed field of characteristic $p>0$, and let $m=p^e>2$ for some integer $e\geq 1$. Let $x$ and $y$ be coordinates on $\AN{2}$ and let $\Ff = \Ll_{\psi, x^{m-1}y}$ be the Artin--Schreier sheaf associated to a nontrivial character $\psi\colon \ZZ/p\ZZ\to \FF_\ell^\times$ and the function $x^{m-1}y$. Let $\pi\colon \AN{2} \to \AN{1}$ be a surjective linear map. The following lemma shows that $R^1 \pi_! \Ff$ is not locally constant. Consequently, the assertion of Proposition~\ref{prop:bertini-nonfierce} is false for sheaves with fierce ramification at infinity.

\begin{lemma} \label{lemma:fierce-is-bad}
	$\Swan_\infty(\Ff|_{\pi^{-1}(0)})\neq \Swan_\infty(\Ff|_{\pi^{-1}(1)})$.
\end{lemma}

\begin{proof}
Say $\pi(x, y) = ax + by$ with $a, b\in k$ not both zero. Suppose first that $b\neq 0$, then $x$ is a coordinate on every fiber of $\pi$, and $y= -\frac{a}{b}x + \frac{1}{b}\pi(x, y)$. If $\pi(x, y) = 0$, then $x^{m-1}y = -\frac{a}{b}x^m$, and hence $\Swan_\infty(\Ff|_{\pi^{-1}(0)}) = 1$ if $a\neq 0$, $0$ if $a=0$. If $\pi(x, y) = 1$, then $x^{m-1}y = -\frac{a}{b}x^m + \frac{1}{b} x^{m-1}$, so  $\Swan_\infty(\Ff|_{\pi^{-1}(1)}) = m-1>1$. It remains to consider the case $b=0$. Then $y$ is a coordinate on every fiber and $x = \frac{1}{a} \pi(x, y)$, so $x^{m-1}y = 0$ if $\pi(x, y) = 0$, and $\Swan_\infty(\Ff|_{\pi^{-1}(0)}) = 0$, and if $\pi(x, y) = 1$ then $x^{m-1} y = \frac{1}{a^{m-1}} y$ and $\Swan_\infty(\Ff|_{\pi^{-1}(1)}) = 1$. In each case we have $\Swan_\infty(\Ff|_{\pi^{-1}(0)}) \neq \Swan_\infty(\Ff|_{\pi^{-1}(1)})$. 
\end{proof}

\begin{corollary}
	In the above situation, $R^q \pi_! \Ff$ is not locally constant for some $q\geq 0$.
\end{corollary}

\begin{proof}
For $t\in \AN 1(k)$, let $\Ff_t = \Ff|_{\pi^{-1}(t)}$. By the Grothendieck--Ogg--Shafarevich formula \eqref{eqn:gos}, we have 
\[ 
	\chi_c(\AN 1, \Ff_t) = \rk(\Ff_t)\cdot \chi_c(\AN 1, \FF_\ell) - \Swan_\infty(\Ff_t) = 1 - \Swan_\infty(\Ff_t).
\]
Thus Lemma~\ref{lemma:fierce-is-bad} implies that the function  $t\mapsto \chi_c(\AN 1, \Ff_t)$ is not constant. On the other hand, 
\[ 
	\chi_c(\AN 1, \Ff_t) = \dim_{\FF_\ell} (R^0 f_! \Ff)_t - \dim_{\FF_\ell} (R^1 f_! \Ff)_t + \dim_{\FF_\ell} (R^2 f_! \Ff)_t, 
\]  
and hence one of the sheaves $R^q f_! \Ff$ is not locally constant.
\end{proof}

\subsection{Complements of hyperplane arrangements}
\label{ex:hyperplane}

Theorem~\ref{thm:affine-kpi1} implies that  every complement of a hyperplane arrangement in $\AN{n}$ is a $K(\pi, 1)$. This is of course false over $\CC$, and this contrast yields examples of interesting arithmetic behavior, the failure of the Lefschetz principle or the question of the existence of a finite \'etale cover killing a given \'etale cohomology class.

\begin{remark}
The question whether certain \emph{complex} complements of hyperplane arrangements are $K(\pi, 1)$ (in the topological sense) has been extensively studied, cf. e.g. \cite{Brieskorn,Deligne,Bessis} or \cite[\S 5.1]{OrlikTerao}. Of course, the fundamental group of the complement of a hyperplane arrangement loses its link to combinatorics (or representation theory) when one passes to positive characteristic. 
\end{remark}

\begin{proposition}
	Let $X = \AA^2_\ZZ \setminus \{ xy(x+y-1)=0\}$. Then
	\begin{enumerate}[(1)]
		\item $\pi_2(X(\CC))\neq 0$ and $\pi_1(X(\CC))\isom \ZZ^3$. In particular, $X(\CC)$ is not a $K(\pi, 1)$ space, and its fundamental group is a good group in the sense of Serre \cite[\S 2.6]{Serre}.
		\item $X_\CC$ is not a $K(\pi, 1)$ scheme. In particular $X_K$ is not a $K(\pi, 1)$ scheme for every field $K$ of characteristic zero.
		\item $X_k$ is a $K(\pi, 1)$ scheme for every field $k$ of characteristic $p>0$.  
	\end{enumerate}
\end{proposition}

\begin{proof}
(1) These statements follow from A.~Hattori's work \cite{Hattori} on the topology complements of generic hyperplane arrangements. See \cite[Example~5.24]{OrlikTerao} for a detailed discussion of this space, whose homotopy type is the same as that of the image $Q$ of the boundary $\partial [0,1]^3$ of the unit cube in $\RR^3$ under the quotient map $\RR^3 \to \RR^3/\ZZ^3$. The inclusion $Q\hookrightarrow \RR^3/\ZZ^3$ induces an isomorphism on fundamental groups, and the surjection $S^2\isom \partial [0,1]^3\to Q$ induces an injection on $\pi_2$. (In \cite[Example~5.24]{OrlikTerao}, the authors claim that this map is an isomorphism on $\pi_2$, which seems to be incorrect.)

(2) 
For a subgroup $\Lambda\subseteq \pi_1(X(\CC))\isom\ZZ^3$ of finite index, let us denote by $X_\Lambda\to X_\CC$ the induced finite \'etale covering. Using the description of (1), $X_\Lambda(\CC)$ can be identified with the image in $\RR^3/\Lambda$ of the union of the boundaries of all unit cubes in $\RR^3$ with vertices in $\ZZ^3$. Let $C$ be the boundary of one of these cubes, then for every finite abelian group $M$ the composition
\[ 
	H^2(X(\CC), M)\to H^2(C, M)\isom H^2(S^2, M)\isom M 
\] 
is easily seen to be surjective. In the commutative diagram
\[ 
\xymatrix{
		H^2(C, M)  & \varinjlim_\Lambda\,\, H^2(X_\Lambda(\CC), M) \ar[l]  & \varinjlim_\Lambda\,\, H^2(X_{\Lambda}, M) \ar[l]_-\sim \\
		& H^2(X(\CC), M)  \ar@{>>}[ul] \ar[u] & H^2(X_\CC, M), \ar[l]_\sim \ar[u]
	}
\]
the slant arrow is surjective, and hence the left vertical arrow is nonzero. Consequently, the top right term is nonzero, which shows that $X_\CC$ is not a $K(\pi, 1)$ by the characterization of Proposition~\ref{kpi1prop} and the fact that $\pi_1(X_\CC)$ is the profinite completion of $\pi_1(X(\CC))$. The case of arbitrary characteristic zero fields follows from \cite[Proposition~3.2(c)]{Achinger}. (In fact, the term $\varinjlim_\Lambda H^2(X_\Lambda, M)$ equals $\Hom(\pi_2(X_\CC), M)$ by \cite[Proposition~6.3]{ArtinMazur}.)

(3) This is a direct consequence of Theorem~\ref{thm:affine-kpi1}.   
\end{proof}

\begin{remark}
An easier example of a scheme presenting such behavior is the `$2$-sphere'
\[ 
	X = \{ x^2 + y^2 + z^2 = 1 \} \subseteq \AA^{3}_{\ZZ[1/2]}.
\]
After base change to a field $k$ containing $\sqrt{-1}$, we can transform the equation to $xy = (z+1)(z-1)$. Since $X$ is affine, $X_k$ is a $K(\pi, 1)$ for every field $k$ of odd characteristic. However, the map $\pi\colon X \to \PP^1$ defined as
\[ 
	\pi(x,y,z) = \begin{cases}
		(x: z-1) & z\neq 1 \\
		(z+1 : y) & z\neq -1
	\end{cases}
\]
makes $X$ into an $\AA^1$-bundle over $\PP^1$. Since $\AA^1_K$ is simply connected when $K$ is algebraically closed of characteristic zero, $X_K$ is simply connected as well. Moreover, the induced map 
\[ 
	\pi^*\colon \ZZ/n\ZZ\isom H^2(\PP^1_{K}, \mu_n)\to H^2(X_{K}, \mu_n)
\]
is an isomorphism, and hence $X_K$ is not a $K(\pi, 1)$.
\end{remark}

\subsection{Fundamental groups of affine spaces} 
\label{ex:differentan}

Let $k$ be an algebraically closed field of characteristic $p>0$. Recall that by the work of Raynaud on Abhyankar's conjecture \cite{Raynaud}, a finite group $G$ arises as a quotient of $\pi_1(\AN{1})$ if and only if $G$ has no nontrivial quotient of order prime to $p$. It follows that $\pi_1(\AN{n})$ has the same property for all $n\geq 1$, and hence the profinite groups $\pi_1(\AN{n})$ and $\pi_1(\AN{1})$ have the same finite quotients. One can ask naively whether $\pi_1(\AN{n})\isom \pi_1(\AN{1})$ as profinite groups (we could deduce this from the previous statement if the groups were topologically finitely generated, cf. \cite[Proposition 15.4]{FriedJarden}). It is easy to see that such an isomorphism cannot be induced by an algebraic morphism $\AN{n}\to \AN{1}$. 

\begin{proposition}
	If $n\neq m$, then $\pi_1(\AN{n})$ and $\pi_1(\AN{m})$ are not isomorphic as profinite groups.
\end{proposition}

\begin{proof}
Theorem~\ref{thm:an-kpi1} implies that the cohomological dimension of $\pi_1(\AN{n})$ equals the largest $q$ for which there exists a locally constant constructible sheaf $\Ff$ on $\AN{n}$ with $H^q(\AN{n}, \Ff)\neq 0$. Thus it suffices to show that the latter equals $n$. This is easy and well-known, but we include a quick proof.

Let $\Ff_1$ be any $\FF_\ell$-sheaf on $\AN{1}$ (for some $\ell\neq p$) with $H^1(\AN{1}, \Ff_1)\neq 0$ (for example, the Artin--Schreier sheaf $\Ll_{\psi, x^m}$ where $m>1$ is an integer prime to $p$ and $\psi\colon  \ZZ/p\ZZ\to \FF_\ell^\times$ is a nontrivial character, cf. \eqref{eqn:gos-no-support} and Proposition~\ref{prop:bryl-nice}).  Let $\Ff_n = \Ff_1\boxtimes \ldots \boxtimes \Ff_1$ ($n$ times). By the K\"unneth formula (cf. \cite[Th. finitude, Corollaire 1.11]{ThFinitude}), $H^n(\AN{n}, \Ff_n) \isom H^1(\AN{1}, \Ff_1)^{\otimes n}\neq 0$. On the other hand, $H^q(\AN{n}, \Ff)$ vanishes for $q>n$ for \emph{all} constructible sheaves $\Ff$, by Artin's theorem on the cohomological dimension of affine schemes \cite[Exp.~XIV, Corollaire~3.2]{Artin}.
\end{proof}

\subsection{Pro-$p$ completion and $p$-Sylow subgroups} \label{ss:pro-p}

Let $\ell$ be a prime. In analogy with the notion of a good group \cite[\S 2.6]{Serre}, let us call a topological group $G$ \emph{$\ell$-good}, if for every continuous representation of the pro-$\ell$ completion $G^{{\rm pro}-\ell}$ on a finite dimensional $\FF_\ell$-vector space $M$, the natural maps
\[ 
	H^*(G^{{\rm pro}-\ell}, M) \to H^*(G, M)
\]
are isomorphisms. 

\begin{proposition}
	Let $(X, \bar x)$ be connected noetherian pointed affine $\FF_p$-scheme. Then the following hold.
	\begin{enumerate}[(1)]
		\item The group $\pi_1(X, \bar x)$ is a $p$-good group,
		\item The $p$-completion $\pi_1(X, \bar x)^{{\rm pro}-p}$ is a free $p$-group of rank equal to $\dim_{\FF_p} \Gamma(X, \Oo_X)/(1-F)$,
		\item Every $p$-Sylow subgroup of $\pi_1(X, \bar x)$ is a free $p$-group. In particular, $\pi_1(X, \bar x)$ is $p$-torsion free. 
	\end{enumerate}
\end{proposition}

\begin{proof}
(1) For brevity, let us call a finite \'etale surjective $f\colon X'\to X$ a \emph{$p$-cover} if it is Galois under the action of a finite $p$-group. Let $M$ be a finite-dimensional $\FF_p$-vector space and let $\rho\colon \pi_1(X, \bar x)\to GL(M)$ be a representation whose image is a $p$-group. Let $\Ff$ be the associated sheaf on $\AN n$. Consider the commutative triangle
\[ 
	\xymatrix{
		H^*(\pi_1(X, \bar x)^{{\rm pro}-p}, M) \ar[dr] \ar[d]  & \\
		H^*(\pi_1(X, \bar x), M) \ar[r] & H^*(X, \Ff). \\
	}
\]
Since the bottom map is an isomorphism by Theorem~\ref{thm:affine-kpi1}, the vertical arrow is an isomorphism if and only if the diagonal one is. To this end, let $\rho_p$ be the natural map of topoi
\[ 
	\rho_p \colon X_\et \to X_{p-f\et} \isom B\pi_1(X, \bar x)^{{\rm pro}-p}.
\]
Here $X_{p-f\et}$ denotes the topos of sheaves on the full subcategory of the \'etale site of $X$ consisting of $p$-covers, with the induced topology. We must show that $R^q \rho_{p*} \Ff = 0$ for $q>0$. As usual, this is automatic for $q=1$. By the definition of the higher direct images, this amounts to showing that for every $p$-cover $X'\to X$, and every class $\zeta \in H^q(X', \Ff)$, there exists a $p$-cover $X''$ of $X$ and a map $f\colon X''\to X'$ over $X$ such that $f^*\zeta = 0 \in H^q(X'', \Ff)$. Since $\Ff$ is trivialized by the $p$-cover corresponding to the image of $\rho$, we can assume $\Ff$ is constant, and hence $\Ff=\FF_p$. In this case $H^q(X, \FF_p) = 0$ for $q>1$ because $X$ is affine, by Artin--Schreier theory. 

(2) By \cite[Theorem~7.7.4]{RibesZaleskii}, to prove that $\pi_1(X, \bar x)^{{\rm pro}-p}$ is free, it suffices to show that $H^2(\pi_1(X, \bar x)^{{\rm pro}-p},\FF_p) = 0$. This follows from (1) and the fact that $H^q(X, \FF_p) = 0$ for $q>1$. The rank can be read off of $\Hom(\pi_1(X, \bar x)^{{\rm pro}-p}, \FF_p) \isom H^1(X, \FF_p) \isom \Gamma(X, \Oo_X)/(1-F)$.

(3) Let $\Pi\subseteq \pi_1(X, \bar x)$ be a $p$-Sylow subgroup. Let $\Pi_\alpha\subseteq \pi_1(X, \bar x)$ be a projective system of open subgroups such that $\Pi = \bigcap_\alpha \Pi_\alpha$, and let $\{X_\alpha\to X\}$ be the corresponding  projective system of finite \'etale coverings of $X$. Since the $X_\alpha$ are also affine $K(\pi, 1)$, we have 
\[ H^2(\Pi, \FF_p) = \varinjlim H^2(\Pi_\alpha, \FF_p) = \varinjlim H^2(X_\alpha, \FF_p) = 0.  \]
Again, we conclude by \cite[Theorem~7.7.4]{RibesZaleskii}.
\end{proof}

\subsection{\texorpdfstring{Relation to $K(\pi, 1)$ pro-$\ell$}{Relation to K(pi, 1) pro-ell}} 
\label{ss:kpi1-pro-ell}

Let $k$ be an infinite field of characteristic $p>0$ and let $\ell\neq p$ a prime. In \cite{FriedlanderKpi1} (see also \cite[Theorem~7.11]{Friedlander}), Friedlander has constructed coverings of smooth schemes $X$ over $k$ by connected open subsets $U$ which are `$K(\pi, 1)$ pro-$\ell$.' This means that the pro-$\ell$ completion of the \'etale homotopy type of $U$ is weakly equivalent to the classifying space of the pro-$\ell$ completion $\pi_1(U,u)^{{\rm pro}-\ell}$. In simpler terms, for every locally constant constructible $\FF_\ell$-sheaf $\Ff$ on $U$ which is `$\ell$-monodromic' (i.e., the image of the associated $\pi_1(U,u)$-representation factors through $\pi_1(U,u)^{{\rm pro}-\ell}$), the natural maps
\[ 
	H^*(\pi_1(U,u)^{{\rm pro}-\ell}, \Ff_u) \to H^*(U, \Ff)
\]
are isomorphisms. Being a $K(\pi, 1)$ does not imply being a $K(\pi, 1)$ pro-$\ell$ \cite[Mise en garde 1.4.5]{MadoreOrgogozo}. A variant of this construction due to Gabber is used in \cite{MadoreOrgogozo} in order to prove computability of $\FF_\ell$-cohomology of schemes of finite type over $k$. Friedlander's construction is a variation of Artin's, and in particular the neighborhoods $U$ he obtains are Artin neighborhoods in the sense defined in the Introduction. 

Since the opens constructed by Friedlander are affine, we deduce that a smooth scheme $X$ over $k$ can be covered by affine open subsets $U$ which are simultaneously $K(\pi, 1)$ and $K(\pi, 1)$ pro-$\ell$ (for a single fixed $\ell$). The two properties imply that for an $\ell$-monodromic $\FF_\ell$-sheaf $\Ff$ on $U$, we have a commutative triangle of isomorphisms
\[ 
	\xymatrix{
		H^*(\pi_1(U, u)^{{\rm pro}-\ell}, \Ff_u) \ar[d] \ar[dr] & \\
		H^*(\pi_1(U, u), \Ff_u) \ar[r]  & H^*(U, \Ff).
	}
\]
This implies that $\pi_1(U, u)$ is $\ell$-good (cf. \S\ref{ss:pro-p}).

\subsection{Anabelian geometry speculations}
\label{ss:anabelian-geometry}

Anabelian geometry started with a series of conjectures formulated by Grothendieck in his letter to Faltings \cite{GrothendieckLetter} (cf. \cite[\S 5]{Szamuely} for a more recent but not completely up to date survey). Its recurring theme is the question whether some class $C$ of connected schemes is `\emph{anabelian},' that is, whether for $X, Y\in C$, the morphism
\[ 
	\Hom(X, Y) \to \Hom_{\rm ext}(\pi_1(X), \pi_1(Y))
\]
is bijective. It only makes sense to consider such classes in characteristic zero, for if $X$ is an $\FF_p$-scheme, then the absolute Frobenius $F_X\colon X\to X$ induces the identity on $\pi_1(X)$. 

To fix this problem, we can introduce an auxiliary category ${\rm Sch}^F$, whose objects are schemes of positive characteristic, and whose morphisms are morphisms of schemes modulo the relation $F_X = \id_X$ ($F_X$ being the absolute Frobenius). More precisely, let ${\rm Sch}_{\rm pos.char.}$ be the full subcategory of the category of schemes ${\rm Sch}$ consisting of schemes $X$ such that $p\cdot \Oo_X$ for some prime $p$ (depending on $X$). Then there is a functor 
\[ 
	\phi\colon {\rm Sch}_{\rm pos.char.} \to {\rm Sch}^F 
\]
such that
\begin{enumerate}[(1)]
	\item $\phi(F_X) = \id_{\phi(X)}$ for all objects $X$,
	\item whenever $\psi\colon {\rm Sch}_{\rm pos.char.} \to \Cc$ is a functor such that $\psi(F_X) = \id_{\psi(X)}$ for all objects $X$, then there exists a unique functor $\bar\psi\colon {\rm Sch}^F \to \Cc$ and an isomorphism $\psi\isom \bar\psi\circ \phi$.
\end{enumerate}
The category ${\rm Sch}^F$ and the functor $\phi$ are uniquely characterized by those properties, up to equivalence. In simple terms, the objects of ${\rm Sch}^F$ are those of ${\rm Sch}_{\rm pos.char.}$, the morphisms $\Hom_{{\rm Sch}^F}(X, Y)$ 
are the quotient of $\Hom(X, Y)$ modulo the relation $f\sim g$ if $F_Y^a\circ f = F_Y^b\circ g$ for some $a, b\geq 0$, and composition is well-defined thanks to the relations $F_Y\circ f = f\circ F_X$.

By the universal property, the functor 
\[ 
	X\mapsto X_\et \colon {\rm Sch}_{\rm pos.char.} \to (\rm Topoi)
\]
descends to ${\rm Sch}^F$, and hence so does the Artin--Mazur homotopy type functor $\Pi$.

\begin{definition} 
A class $C$ of objects of ${\rm Sch}_{\rm pos.char.}$ is \emph{anabelian} if for $X, Y\in C$, the morphism
\begin{equation} \label{eqn:anabelian-map} 
	\Hom_{{\rm Sch}^F}(X, Y)\to \Hom_{\rm pro-H}(\Pi(X), \Pi(Y))
\end{equation}
is a bijection.
\end{definition}

We put $\Pi(X)$ here instead of $\pi_1(X)$ because it behaves better in the absence of basepoints (cf. \cite[\S 2.2]{SchmidtStix}), and because the schemes we are interested in are $K(\pi, 1)$ anyway. 

The following question seems very natural in the light of our results.

\begin{question}
	Let $C$ be the class of integral normal affine schemes $X=\Spec R$ with $p\cdot \Oo_X = 0$ for some prime $p$ depending on $X$, such that if $R^{p^\infty} = \bigcap_{e\geq 0} F_R^e(R)$, then
	\begin{enumerate}[(1)]
		\item the perfect ring $R^{p^\infty}$ is a field,
		\item $R$ is finitely generated over $R^{p^\infty}$,
		\item $\dim X > 0$. 
	\end{enumerate}
	Is $C$ an anabelian class?
\end{question}

Note the contrast with Grothendieck's conjectures, where the fields in question are finitely generated. In fact, it makes sense to ask the above question under the additional restriction that $R^{p^\infty}$ is an algebraically closed field.

It is not difficult to show, using the Artin--Schreier isomorphism $\Hom(\pi_1(X), \FF_p)\isom H^1(X, \FF_p) \isom R/(1-F)$, that the map \eqref{eqn:anabelian-map} is always injective. In other words, if $X$ and $Y$ are in the class $C$, and $f, g\colon X\to Y$ are two morphisms inducing the same map $H^1(Y, \FF_p)\to H^1(X, \FF_p)$, then $f=g$. The details will appear elsewhere.

Surjectivity of \eqref{eqn:anabelian-map} seems beyond reach at the present moment. A seemingly more tractable special case is the following:

\begin{question}
 	Let $k$ be a perfect field of characteristic $p$. Does $\pi_1(\AN 1)$ determine $k$? 
\end{question} 

\appendix

\section{The rank one case of the Bertini theorem} 
\label{s:appendix}

We present a second proof of Theorem~\ref{thm:bertini} in the case when the sheaf $\Ff$ has rank one. In this case, we are able to explicitly find the desired automorphism $\phi$. Interestingly, we do not have to assume that $k$ is infinite.

\subsection{Ramification of Artin--Schreier--Witt sheaves}

For $r\geq 0$, let $\WW_{r+1}$ denote the group scheme of Witt vectors of length $r+1$ over $k$, and let $F\colon \WW_{r+1}\to \WW_{r+1}$ be the Frobenius homomorphism. As a variety, $\WW_{r+1}\isom \AN{r+1}$, and $F(x_0, \ldots, x_{r}) = (x_0^p, \ldots, x_{r}^p)$. We have a short exact (for the \'etale topology) sequence of $k$-group schemes:
\begin{equation} \label{eqn:witt-ses} 
	0\to \ZZ/p^{r+1}\ZZ\to \WW_{r+1} \xto{1-F} \WW_{r+1} \to 0,
\end{equation}
yielding for any $k$-scheme $X$ a boundary map
\begin{equation} \label{eqn:delta}
	\delta \colon \WW_{r+1}(X) \to H^1(X, \ZZ/p^{r+1}\ZZ).
\end{equation}
If $H^1(X, \Oo_X)=0$ (e.g. if $X$ is affine), then $\delta$ induces an isomorphism
\[ 
	\WW_{r+1}(X)/(1-F)\WW_{r+1}(X) \isomlong H^1(X, \ZZ/p^{r+1}\ZZ). 
\]

\begin{definition} \label{def:asw}
	If $\ell\neq p$ is a prime and $\psi\colon \ZZ/p^{r+1}\ZZ \to \bar\FF_\ell^\times$ is an injective character, we denote by $\Ll_\psi$ the rank one $\bar\FF_\ell$-sheaf on $\WW_{r+1}$ associated to the $\ZZ/p^{r+1}\ZZ$-torsor given by \eqref{eqn:witt-ses} and the character $\psi$. 

	If $X$ is a $k$-scheme and $f\colon X\to \WW_{r+1}$ is a map, we denote by $\Ll_{\psi,f}$ the pull-back $f^* \Ll_\psi$; we call it the \emph{Artin--Schreier--Witt sheaf} associated to $f$ and $\psi$. 
\end{definition}

It is clear that if $f-h \in (1-F)\WW_{r+1}(X)$, then $\Ll_{\psi,f}\isom \Ll_{\psi,h}$. Our first goal in this section is to study the ramification of $\Ll_{\psi,f}$ along a divisor, in particular its Swan conductor. 

To this end, suppose that $X=\eta=\Spec K$ as in Section~\ref{s:wild}. In this context, and in the perfect residue field case ($\kappa=k$), the Swan conductor of an Artin--Schreier--Witt sheaf over $X$ has been computed by Brylinski \cite{Brylinski} and Kato \cite{KatoSwan}. To explain the result, we endow the groups $\WW_{r+1}(K)$ and $H^1(K, \ZZ/p^{r+1}\ZZ)$ with increasing filtrations as follows:
\begin{align} 
	\label{eqn:filww}
	\fil_n \WW_{r+1}(K) &= \left\{ f=(f_0, \ldots, f_r)\in \WW_{r+1}(K)\,:\, -p^{r-i}\nu(f_i) \leq n \right\},\\
	\label{eqn:filh1}
	\fil_n H^1(X, \ZZ/p^{r+1}\ZZ) &= \delta\left(\fil_n \WW_{r+1}(K)\right), 
\end{align}
where $\nu$ is the valuation on $K$. 

\begin{theorem}[{\cite[Corollary to Theorem 1]{Brylinski}}] \label{thm:bryl}
	Suppose that the residue field of $K$ is $k$. Let $f=(f_0, \ldots, f_{r})\in \WW_{r+1}(K)$, $\psi\colon \ZZ/p^{r+1}\ZZ \to \bar\FF_\ell^\times$ an injective character, defining the $\bar\FF_\ell$-sheaf $\Ll_{\psi, f}$ on $\Spec K$.  Then
	\begin{align*} 
		\Swan(\Ll_{\psi, f}) &= \min \left(\{0\}\cup\left\{ n \,:\, \delta(f)\in \fil_n H^1(K, \ZZ/p^{r+1}\ZZ)\right\}\right) \\
			&=   \min_{h} \max \{ 0, -p^{r} \nu(h_0), -p^{r} \nu(h_1), \ldots, -\nu(h_{r}) \},    
	\end{align*}
	where $h$ ranges over all $h\in \WW_{r+1}(K)$ such that $f-h\in (F-1)\WW_{r+1}(K)$.
\end{theorem}

The $\min_h$ in the above theorem might be difficult to compute in concrete situations. The following proposition provides a simple condition for that minimum to be attained at $h=f$.

\begin{proposition} \label{prop:bryl-nice}
	In the situation of Theorem~\ref{thm:bryl}, suppose that 
	\begin{enumerate}[(1)]
		\item the function $i\mapsto -p^{r-i}\nu(f_i)$ has a unique maximum $i_0$,
		\item the number $\nu(f_{i_0})$ is negative and prime to $p$.
	\end{enumerate}
	Then $\Swan(\Ll_{\psi, f}) = -p^{r-i_0}\nu(f_{i_0})$.
\end{proposition}

Following Abbes and Saito \cite[\S 10]{AbbesSaitoMicrolocale}, we endow $\Omega^1_{K/k}$ with the increasing filtration 
\begin{equation} \label{eqn:filomega} 
	\fil_n \Omega^1_{K/k} = {\rm image}(\Omega^1_{\Oo_K/k}(\log)\otimes_{\Oo_K} \mm_{\Oo_K}^{-n} \to \Omega^1_{K/k} ). 
\end{equation}
Here $\Omega^1_{\Oo_K/k}(\log)$ is the module of log differentials (cf. \cite[\S 5.4]{AbbesSaitoMicrolocale}). The map
\[ 
	F^r d \colon \WW_{r+1}(K) \to \Omega^1_{K/k}, \quad
	f \mapsto \sum_{i=0}^r f_i^{p^{r-1}-1}df_i
 \]
is additive and compatible with the filtrations (the notation $F^r d$ comes from the de~Rham--Witt complex). Moreover, upon passing to the associated graded quotients (${\rm Gr}_n = \fil_n/\fil_{n-1}$), $F^r d$ factors uniquely through $\delta$ \eqref{eqn:delta}, that is, for every $n$ there exists a unique map $\psi_n$ making the following triangle commute
\[ 
	\xymatrix{
		\Gr_n \WW_{r+1}(K) \ar[r]^{F^r d} \ar[d]_{\Gr_n} & \Gr_n \Omega^1_{K/k}  \\
		\Gr_n H^1(K, \ZZ/p^{r+1}\ZZ).\ar[ur]_{\psi_n}
	}
\]
Furthermore, $\psi_n$ is injective. 

\begin{lemma} \label{lemma:abc}
	Let $f=(f_0, \ldots,f_r)\in \WW_{r+1}(K)$ and 
	\begin{align*}
		a &= \min\{ n\,:\, f\in \fil_n \WW_{r+1}(K) \},  & \text{(cf. \eqref{eqn:filww})}\\
		b &= \min\{ n\,:\, \delta(f)\in \fil_n H^1(K, \ZZ/p^{r+1}\ZZ) \},& \text{(cf. \eqref{eqn:filh1})}\\
		c &= \min\{ n\,:\, F^r d (f)\in \fil_n \Omega^1_{K/k}\} & \text{(cf. \eqref{eqn:filomega})}.
	\end{align*}
	Then $a=c$ if and only if $a=b$.
\end{lemma}

\begin{proof}
This follows from the existence and injectivity of $\psi_n$. 
\end{proof}

\begin{proof}[Proof of Proposition~\ref{prop:bryl-nice}]
	Let $a,b,c$ be as in Lemma~\ref{lemma:abc}. Then $a = -p^{r-i_0}\nu(f_{i_0})$, and Theorem~\ref{thm:bryl} states that $b=\Swan(\Ll_{\psi, f})$. Moreover, 
	\[ 
		F^rd(f) = \left( \alpha\nu(f_{i_0})t^{p^{r-i_0}\nu(f_{i_0})} + \text{higher order terms}\right) d\log t
	\]
	with $\alpha, \nu(f_{i_0})\in k^*$. Thus $c=-p^{r-i_0}\nu(f_{i_0})=a$, so $a=b$.
\end{proof}

\subsection{Artin--Schreier--Witt sheaves on $\AN n$}

\begin{lemma} \label{lemma:uniqmax}
	Let $p$ be a prime number and let $S$ be a non-empty finite subset of $\NN^{n+1}$. Let $\leq$ be the order on $S$ such that $(b_0, b_1, \ldots, b_n)\leq (b'_0, b'_1, \ldots, b'_n)$ if and only if $(b_0+b_1+\ldots+b_n, b_0, b_1, \ldots, b_n)\leq (b'_0+b'_1+\ldots+b'_n, b'_0, b'_1, \ldots, b'_n)$ in the lexicographical order. This is a linear order on $S$; let $x=(b^*_0, b_1^*, \ldots, b_n^*)$ be the largest element of $S$. There exist integers $d_1>\ldots>d_n>0$ such that if we set $w(b_0, b_1, \ldots, b_n) = b_0 + d_1 b_1 + \ldots + d_n b_n$ then
	\begin{enumerate}
  		\item $x$ is the unique maximum of $w$ on $S$,
  		\item $w(x)/{\rm gcd}(b_0^*, b_1^*, \ldots, b_n^*)$ is not divisible by $p$.
 	\end{enumerate}  
\end{lemma}

\begin{proof}
Let $A$ be the diameter of $S$ for the $\ell^\infty$-norm, so if $(b_0, b_1, \ldots, b_n)$, $(b'_0, b'_1, \ldots, b'_n) \in S$ then $A \geq |b_i-b'_i|$. It is easy to see that $w\colon S\to \ZZ$, $w(b_0, b_1, \ldots, b_n) = b_0 + d_1 b_1 + \ldots + d_n b_n$, will be injective as long as $d_{i+1} > A(1 + d_1 + \ldots + d_i)$ for all $i=1, \ldots, n-1$. Moreover, the condition that $d_1>\ldots>d_n>0$ ensures that $w\colon S\to \ZZ$ is non-decreasing. This implies that for every tuple of integers $e_1, \ldots, e_n$, we can find a tuple $d_1> \ldots > d_n > 0$ satisfying (1) and such that $d_i \equiv e_i$ modulo $p$. We deduce that we can make condition (2) fulfilled as well.
\end{proof}

\begin{lemma} \label{lemma:finddi}
	Let $f_0, \ldots, f_{r} \in k[x_0, \ldots, x_n]$ be polynomials without non-constant $p$-power monomials, such that not all $f_i$ are constant. Then there exist $d_1 >  \ldots > d_n>0$ such that if 
	\[ 
		f'_i = f_i(x_0, x_1 + x_0^{d_1}, \ldots, x_{n} + x_0^{d_{n}}), 
	\]
	then 
	\begin{enumerate}[(1)]
		\item $i\mapsto \deg(f'_i)p^{r-i}$ has a unique minimum $\deg(f'_{i_0})p^{r-i_0}$, and 
		\item $p$ does not divide $\deg(f'_{i_0})$.
	\end{enumerate}
\end{lemma}

\begin{proof}
Apply Lemma~\ref{lemma:uniqmax} to $S = \bigcup {\rm supp}(f_i)\cdot p^{r-i}$. By assumption, this union is disjoint.
\end{proof}

\begin{lemma} \label{lemma:noppower}
	Let $A=k[x_0, \ldots, x_n]$, $f\in \WW_{r+1}(A)$. Then there exists a $g\in \WW_{r+1}(A)$ such that if $f + (1-F)g = (h_0, \ldots, h_{r})$ then the $h_i$ have no nonzero $p$-th power monomials.
\end{lemma}

\begin{proof}
The proof is by induction on $r\geq 0$. Elementary combinatorics shows that for $u\in A$, there exists a $v(u) \in A$ (unique up to adding an element of $\FF_p$) such that $w(u) = u + v(u) - v(u)^p$ has no non-zero $p$-th power monomials (note that we need $k$ separably closed to deal with the constant term). This solves the base case $r=0$. For the induction step, we take $f= (f_0, \ldots, f_{r})$ with $r>0$, and set 
\[ f' = f + (1-F)(v(f_0), 0, \ldots, 0). \]
Clearly $f=f'$ mod $(1-F)\WW_{r+1}(A)$. Moreover, $f'_0 = f_0 + v(f_0) - v(f_0)^p = w(f'_0)$ has no $p$-th power monomials. Set
\[ f''= (f'_1, \ldots, f'_{r}) \in \WW_{r}(A), \]
so that $f' = (f'_0, 0, \ldots, 0) + Vf''$. By the induction assumption, there exists a $g''$ such that the entries of $h'' = f'' + (1-F)g''$ have no $p$-th power monomials. Finally, let 
\[ g = (v(f_0), 0, \ldots, 0) + (0, g''_0, \ldots, g''_{r-1}) = (v(f_0), 0, \ldots, 0) + Vg''. \]
Then 
\begin{align*} 
	h = f + (1-F)g &= f + (1-F)(v(f_0), 0, \ldots, 0) + (1-F)Vg'' \\
	&= f' + (1-F)Vg'' = (f'_0, 0, \ldots, 0) + Vf'' + V(1-F)g'' \\
	&= (w(f'_0), h''_0, \ldots, h''_{r-1}),
\end{align*}
and we see that the entries of $h$ have no $p$-th power monomials, as desired.
\end{proof}

\begin{proposition}
	The assertion of Theorem~\ref{thm:bertini} holds if $\Ff$ is a rank one $\bar\FF_\ell$-sheaf.
\end{proposition}

\begin{proof}
As $\pi_1(\AN{n+1})$ has no prime-to-$p$ quotients, while $\bar\FF_\ell^\times$ is cyclic, the representation $\pi_1(\AN{n+1})\to \bar\FF_\ell^\times$ corresponding to $\Ff$ factors through an injective character $\psi\colon \ZZ/p^{r+1}\ZZ\to \bar\FF_\ell^\times$ for some $r\geq 0$. Thus there exists an $f=(f_0, \ldots, f_{r}) \in \WW_{r+1}(k[z, x_1, \ldots, x_n])$ such that $\Ff \isom \Ll_{\psi, f}$. We can assume that not all of the $f_i$ are constant, otherwise $\Ff$ is constant and there is nothing to prove. By Lemma~\ref{lemma:noppower}, we can moreover assume that the $f_i$ have no nonzero $p$-th power monomials, as $\Ll_{\psi, f}\isom \Ll_{\psi, f+(1-F)g}$ for any $g\in \WW_{r+1}(k[z, x_1, \ldots, x_n])$. By Lemma~\ref{lemma:finddi}, there exist $d_1 >  \ldots > d_n>0$ such that if $f'_i = f_i(z, x_1 + z^{d_1}, \ldots, x_{n} + z^{d_{n}})$, then $i\mapsto \deg(f'_i)p^{r-i}$ has a unique minimum $\deg(f'_{i_0})p^{r-i_0}$, and $p$ does not divide $\deg(f'_{i_0})$. But by Proposition~\ref{prop:bryl-nice} this means that $\Swan_\infty(\Ff_{\pi^{-1}(y)}) = \deg(f'_{i_0})p^{r-i_0}$ is independent of $y\in \AN{n}$, where $\pi = {\rm pr}\circ \phi$, 
\[
	\phi(z, x_1, \ldots, x_n) = (z, x_1 + z^{d_1}, \ldots, x_n + z^{d_n})
\] 
and ${\rm pr}(z, x_1, \ldots, x_n) = (x_1, \ldots, x_n)$. 
\end{proof}

\begin{remark}
We expect that the sheaf $\phi^* \Ff$ in the proof above is non-fierce at infinity.
\end{remark}

\begin{remark}
Our calculations (especially Lemma~\ref{lemma:uniqmax}) resemble those of Barrientos \cite[\S 5]{Barrientos} in a similar context of restricting rank one sheaves to curves.
\end{remark}

\bibliographystyle{amsalpha} 
\bibliography{article}

\end{document}